\documentclass[11pt]{article}

\usepackage[usenames,dvipsnames]{pstricks}
\usepackage{epsfig}
\usepackage{pst-grad} 
\usepackage[space]{grffile} 
\usepackage{etoolbox} 
\makeatletter 
\patchcmd\Gread@eps{\@inputcheck#1 }{\@inputcheck"#1"\relax}{}{}
\makeatother

\usepackage{graphicx}
\usepackage{amsmath}
\usepackage{psfrag}
\usepackage{color,soul}
\usepackage{siunitx}
\usepackage{float}
\usepackage{overpic}
\usepackage{tikz}
\usepackage{multirow}
\usepackage{bm}
\usepackage{textgreek}
\usepackage{url}
\usepackage{amssymb}
\usepackage{wasysym}
\usepackage{xfrac}
\usepackage{setspace}
\usepackage{gensymb}
\usepackage{mathtools}
\usepackage{arydshln}
\usepackage{algorithm}
\usepackage{algpseudocode}

\usepackage[noblocks]{authblk}
\usepackage[margin=1in]{geometry}
\usepackage[title]{appendix}

\newcommand\affiliation[1]{\gdef\@affiliation{\let\aff\aff@inst#1}}
\gdef\@affiliation{}
\def\email#1{Email address for correspondence: #1}
\def\aff#1{\ignorespaces\textsuperscript{#1}}
\def\corresp#1{\unskip\thanks{#1}}
\numberwithin{equation}{section}

\renewenvironment{abstract}
{\begin{quote}
\noindent \rule{\linewidth}{.5pt}\par{\bfseries \abstractname.}}
{\medskip\noindent \rule{\linewidth}{.5pt}
\end{quote}
}

  \DeclareTextFontCommand\textsfi{\usefont{OT1}{cmss}{m}{sl}}
  \DeclareMathAlphabet\mathsfi            {OT1}{cmss}{m}{sl}
  \DeclareTextFontCommand\textsfb{\usefont{OT1}{cmss}{bx}{n}}
  \DeclareMathAlphabet\mathsfb            {OT1}{cmss}{bx}{n}
  \DeclareTextFontCommand\textsfbi{\usefont{OT1}{cmss}{m}{sl}}
  \DeclareMathAlphabet\mathsfbi            {OT1}{cmss}{m}{sl}
  \IfFileExists{t1phv.fd}
  {\DeclareTextFontCommand\textsfbi{\usefont{T1}{phv}{b}{it}}
  \DeclareMathAlphabet\mathsfbi            {T1}{phv}{b}{it}
  }{}
  \IfFileExists{ot1phv.fd}
  {\DeclareTextFontCommand\textsfbi{\usefont{OT1}{phv}{b}{it}}
  \DeclareMathAlphabet\mathsfbi            {OT1}{phv}{b}{it}
  }{}

\newcommand{\td}[1]{ {\bf #1} }			
\newcommand{\tdh}[1]{ \hat{\td{#1}}}
\newcommand{\vd}[1]{ {\boldsymbol #1} } 
\newcommand{\vdh}[1]{ \hat{\vd{#1}}}    
\newcommand{\vdt}[1]{ \tilde{\vd{#1}}}    
\newcommand{\SPOD}[1]{\td{\Psi}_{#1}}   
\newcommand{\set}[1]{\mathcal{#1}}

\makeatletter
\newcommand{\lambdabar}{{\mathchoice
  {\smash@bar\textfont\displaystyle{0.25}{1.2}\lambda}
  {\smash@bar\textfont\textstyle{0.25}{1.2}\lambda}
  {\smash@bar\scriptfont\scriptstyle{0.25}{1.2}\lambda}
  {\smash@bar\scriptscriptfont\scriptscriptstyle{0.25}{1.2}\lambda}
}}
\newcommand{\smash@bar}[4]{%
  \smash{\rlap{\raisebox{-#3\fontdimen5#10}{$\m@th#2\mkern#4mu\mathchar'26$}}}%
}
\makeatother

\DeclareSymbolFont{matha}{OML}{txmi}{m}{it}
\DeclareMathOperator*{\argmax}{arg\,max}

\usepackage{hyperref}

\definecolor{darkblue}{rgb}{0,0,0.80}

\hypersetup{
    colorlinks=true,
    linkcolor={darkblue},
    citecolor={darkblue},
    filecolor=black,
    urlcolor={darkblue},
    plainpages=false,
}

\providecommand{\keywords}[1]
{

  \textbf{Key words: } #1
}
\def\showfigs{1} 

\title{{\bf Nonlinear space-time model reduction in the frequency domain}}

\author[1]{\bf Peter Frame\corresp{\email{pframe@umich.edu}}}
\author[1]{\bf Aaron Towne}

\affil[1]{\normalsize Department of Mechanical Engineering, University of Michigan, Ann Arbor, MI, USA  }

\date{}
\begin{document}
\maketitle

\begin{abstract}
We propose a space-time reduced-order model (ROM) for nonlinear dynamical systems, building upon previous work on linear systems. Whereas most ROMs are space-only in that they reduce only the spatial dimension of the state, the proposed method leverages an efficient encoding of the entire trajectory of the state on the time interval $[0,T]$, enabling significant additional reduction. Trajectories are encoded using SPOD modes, a spatial basis at each temporal frequency tailored to the structures that appear at that frequency. These modes have a number of properties that make them an ideal choice for space-time model reduction, including separability and near-optimality for long trajectories. We derive a system of algebraic equations involving the SPOD coefficients, forcing, and initial condition by projecting an implicit solution of the governing equations onto the set of SPOD modes in a space-time inner product. We therefore refer to the method as spectral solution operator projection (SSOP). The online phase of SSOP comprises solving this system for the SPOD coefficients, given the initial condition and forcing. We find that SSOP gives two orders of magnitude lower error than POD-Galerkin projection at the same number of modes and CPU time across a suite of tests, including ones that use out-of-sample forcings and affine parameter variation. In fact, the method is substantially more accurate even than the projection of the solution onto the POD modes, which is a lower bound for the error of any method based on a linear space-only encoding of the state.
\\
\end{abstract}


\section{Introduction}
\label{Sec:Intro}
The great majority of model reduction methods use the following two-step approach: \textit{(i)} find an efficient encoding of the state of the dynamical system, replacing the high-dimensional state space with a low-dimensional one, then \textit{(ii)} derive a dynamical system for the encoding. There are many methods for both of these steps. Proper orthogonal decomposition (POD) modes \cite{Lumley67,Sirovich87} are perhaps the most common choice for \textit{(i)}, but other choices include balanced truncation modes \cite{Moore81,Willcox02,Rowley05} and nonlinear encodings like those based on autoencoders \cite{Choi19}. These methods all make use of correlations between the points that make up the state. Common choices for \textit{(ii)} include methods like Galerkin \cite{Aubry88,Noack03} and Petrov-Galerkin \cite{Carlberg10} projections, which derive the reduced-order model (ROM) equations from the full-order model (FOM) equations, and methods like operator inference \cite{Peherstorfer16,Kramer24}, which obtain equations for the evolution of the encoding directly from data. The former methods are referred to as intrusive, and the latter as non-intrusive. We refer to methods that use this two-step paradigm as space-only model reduction methods --- only the spatial dimension of the system is reduced.
\\

Far less common are space-time model reduction methods \cite{Yano14,Choi19,Choi21,Choi21b,Parish21,Hoang22}, which leverage the spatiotemporal correlations present in problems of interest, rather than solely the spatial ones. Step \textit{(i)} is now to find an efficient encoding for the entire trajectory of the state over some time interval of interest. With the entire evolution of the state represented in the encoding, the degrees of freedom that represent the solution no longer evolve in time. Step \textit{(ii)} becomes solving for these static degrees of freedom either by solving a system of simultaneous equations or an optimization. The potential advantage of the space-time paradigm over the space-only one is that solutions to problems of interest exhibit spatiotemporal coherent structures, which offer an extremely efficient means of describing the dynamics. In a space-only ROM, one is nonetheless solving for an encoding of the trajectory -- the space-only encoding of the state at all time steps may be viewed as such -- but this trajectory encoding is far from optimal. Indeed, the space-only POD coefficients, for example, are highly correlated from one time step to the next, so they provide a redundant (i.e., suboptimal) encoding of the entire trajectory. Thus, with the same total degrees of freedom, the space-time framework offers the potential for substantial additional accuracy by leveraging spatiotemporal correlations. The key questions are the following. (1) Can one solve for the space-time encoding coefficients accurately? POD-Galerkin projection (POD-G), e.g., does not exactly recover the projection of the solution onto the POD modes due to the influence of unclosed terms. To what extent does this persist for a space-time approach? (2) Can the system of equations in the space-time ROM be solved in similar time to the time integration required in the space-only ROM? In other words, does the space-time approach take the same (or less) CPU time as the space-only approach for the same number of degrees of freedom? If the answer to both of these is `yes', then the space-time ROM will be much more accurate than a space-only ROM at the same CPU time.
\\

The spatiotemporal basis vectors used in many space-time model reduction methods are separable, i.e., they are formed by taking the Kronecker product of a spatial vector and a temporal one, and the spatial basis vectors used are often the standard POD modes \cite{Choi19,Choi21,Choi21b,Hoang22}. In most cases, the time dependence for a given POD mode is tailored to that mode, which is accomplished by taking the SVD of a data matrix containing as its columns different temporal evolutions of the associated POD coefficient. In other cases \cite{Choi19}, the same temporal basis is used for each POD mode. In either case, the temporal basis is obtained empirically. Refs. \cite{Choi19,Choi21,Choi21b} minimize a space-time residual over the basis coefficients, and Ref. \cite{Hoang22} trains a machine learning model that outputs the basis coefficients.
\\

Our approach is to use spectral proper orthogonal decomposition (SPOD) modes as the space-time basis. At each temporal frequency, there is a spatial basis of modes --- the SPOD modes --- that are optimal for describing the structures that appear at that temporal frequency. The modes at frequency $\omega_k$ are associated, by definition, with the time dependence $e^{i\omega_kt}$. Therefore, the SPOD modes can be used as a (separable) space-time basis by adding these modes with their associated time dependence together with constant coefficients. A number of properties of the SPOD basis make them attractive for model reduction. First, for statistically stationary systems, SPOD modes are nearly the optimal linear space-time encoding for long time intervals. More precisely, for statistically stationary systems, SPOD modes approach space-time POD modes, the optimal linear encoding of trajectories, as the time interval of the encoding goes to infinity \cite{Lumley67,Towne2018spectral}. Second, various properties of the analytic time dependence of the modes are useful; one of particular importance is that the Fourier time dependence leads to sparser coupling via the nonlinearity than one would achieve with a general space-time basis. Finally, the SPOD basis is separable, i.e., each mode is a function of space multiplied by a function of time. This means that relative to other space-time bases such as space-time POD, the SPOD modes are easier to converge from data, consume less memory, and inner products are far cheaper to evaluate.
\\

In our previous work \cite{Frame25}, we proposed \textit{Spectral Solution Operator Projection} (SSOP), a space-time ROM that returns SPOD coefficients that represent the solution to a linear time-invariant system on a predefined time interval, given the initial condition and forcing. We showed that this method achieved orders-of-magnitude lower error than POD-G, balanced truncation \cite{Moore81}, and even the projection of the FOM solution onto POD modes. Such accuracy was achieved at comparable (and usually slightly lower) CPU time than these benchmarks. SSOP is so named because it uses the fundamental solution to LTI systems, which can be thought of as a solution operator acting on the initial condition and forcing, and applies a space-time projection onto spectral POD modes to this solution operator. This results in an equation, at every frequency, giving the SPOD coefficients at that frequency in terms of the forcing and initial condition.
\\

The goal of this paper is to generalize SSOP to nonlinear systems. Three new aspects must be addressed. First, the operator must be adjusted to account for the nonlinearity. To do this, we treat the nonlinear term as an additional forcing on the system and use the fundamental solution to the linear system with this additional forcing. This is consistent with an increasingly popular perspective in the fluid dynamics community on the role of nonlinearity \cite{McKeon10,McKeon17}. Crucially, this additional forcing depends on the solution itself, and this results in a system of nonlinear equations that must be solved for the SPOD coefficients. Second, the influence of nonlinearity at each frequency must be approximated as a function of the retained SPOD coefficients. We use two approaches, depending on the nature of the nonlinearity. The first is a DEIM-based hyper-reduction approximation of the nonlinearity, and the second, which is only applicable for quadratic nonlinearities, uses a subset of the triadic interactions between the SPOD modes. Third, the system of nonlinear equations that results from the addition of nonlinearity must be efficiently solved. We use a simple fixed-point iteration and find that for most cases, this converges in $\sim 10$ iterations. For cases where this does not converge, we use a pseudo-time-stepping method, which is somewhat slower but more stable.
\\

We test the method on two Ginzburg-Landau systems, one with the standard cubic nonlinearity and the other with a Burgers'-type (quadratic) nonlinearity used to test the triadic-interaction handling of the nonlinearity. We demonstrate good results for both: at similar cost to POD-G, SSOP results in orders-of-magnitude lower error than POD-G and even the lower bound for space-only methods given by the projection of the solution onto the POD modes. We test the method both on statistically stationary forcings, where the SPOD modes provide the asymptotically optimal encoding, and on non-stochastic forcings, observing similar results for both. We also test the ability of the method to handle affine parameter dependence --- we use SPOD modes from one Ginzburg-Landau system to build the ROM, then test it on a parameter sweep of Ginzburg-Landau systems.
\\

In addition to our previous work \cite{Frame25}, several methods have been proposed that use SPOD modes for model reduction. Refs. \cite{Lin19, Towne21} developed Galerkin and Petrov-Galekin SPOD-based ROMs, respectively, for linear systems. These models are only capable of handling systems in which the solution is known to be periodic. Recently, Ref. \cite{Li25} proposed a nonlinear SPOD-based model, somewhat similar to the Galerkin models in Refs. \cite{Lin19,Towne21}, but which was generalized to handle nonlinearities. Ref. \cite{Li25} then applied this model to find periodic solutions in a canonical incompressible viscous flow. This latter model can be seen as a spatially reduced version of the harmonic balance method \cite{Hall02}, a method which derives a relation between the various frequencies of a temporally periodic solution to the Navier-Stokes equations. The harmonic balance method has been applied in turbomachinery problems \cite{Hall02,Hall13}, where the solution is assumed to be periodic with the fundamental frequency given by the blade-passing frequency, and has been used to study boundary layer transition \cite{rigas21}. In both, just a few harmonics are retained. All the methods above require a periodic solution, something that might be natural to expect of methods based on the Fourier transform. However, SPOD modes are effective at representing both periodic and aperiodic trajectories, and since the latter is far more common, there is a greater need for reduced-order models for that case. The proposed method, and its antecedent \cite{Frame25}, can handle both periodic and aperiodic problems, owing to the solution operator projection approach.
\\

The remainder of this paper is organized as follows. In Section~\ref{sec:preliminaries}, we describe trajectory representation, SPOD modes, and linear SSOP. We derive the general nonlinear SSOP in Section~\ref{sec:method}. We report the results of the application of the method to the Ginzburg-Landau systems in Section~\ref{sec:results}, and finally conclude the paper in Section~\ref{sec:conclusions}.
\\

\section{Spectral POD and linear spectral solution operator projection} \label{sec:preliminaries}
In this section, we discuss the preliminaries that lay the foundation for the proposed method. First, we discuss the task of trajectory encoding. Then, we review the relevant properties of spectral POD \cite{Towne2018spectral}, the most important of which is the way SPOD modes can be used to represent a trajectory and the efficiency of this representation \cite{Frame23}. Finally, we review linear SSOP \cite{Frame25}, which is the antecedent to this work. 
\\

\subsection{Trajectory representation}
A common task in model reduction is formulating an accurate scheme for encoding the state of a dynamical system $\vd{q} \in \mathbb{C}^{N_x}$. This entails specifying an encoder $\vd{e}^{x}: \mathbb{C}^{N_x} \to \mathbb{C}^r$ and a decoder $\vd{d}^x: \mathbb{C}^{r} \to \mathbb{C}^{N_x}$, where, ideally, $r \ll N_x$ and $\vd{q} \approx \vd{d}^x(\vd{e}^x(\vd{q}))$ for some appropriate ensemble of states $\vd{q}$. The superscript indicates that the spatial dimension is reduced. It is well known that, in a certain precise sense, the first $r$ proper orthogonal decomposition modes can be used to formulate the optimal linear encoding scheme. Specifically, the encoder $\vd{e}^{x}(\vd{q}) = \td{\Phi}^*\td{W} \vd{q}$ and decoder $\vd{d}^{x}(\vd{a}) = \td{\Phi}\vd{a}$ minimize $\mathbb{E}[\| \vd{q} -  \vd{d}^x(\vd{e}^x(\vd{q})) \|_{x}^2]$ over all linear encoder/decoder pairs \cite{Oja82}, where $\td{\Phi} \in \mathbb{C}^{N_x \times r}$ is the matrix of the first $r$ POD modes, and $\td{W}$ is a (positive definite) weight matrix. The norm is induced by the inner product $\| \vd{q} \|_{x} = \sqrt{\langle \vd{q}, \vd{q} \rangle_{x}}$, the inner product is defined using the weight matrix as $\langle \vd{q}_1, \vd{q}_2 \rangle_{x} = \vd{q}_2^* \td{W} \vd{q}_1$, and $\mathbb{E}[ \cdot ]$ is the expectation operator over the ensemble that defines the POD modes. The `$x$' subscript is used to denote a space-only norm or inner product, as opposed to a space-time one.
\\

A related, though less familiar, task is formulating an encoding/decoding scheme for trajectories of the state of the dynamical system on some time interval $[0,T]$. We take the trajectory on the interval to be a vector of its values at times $0, \Delta t, \dots, (N_\omega-1)\Delta t$, where $N_\omega \Delta t = T$. That is, the trajectory is a vector $\vd{q}_\set{J} = [\vd{q}_0^T, \vd{q}_1^T, \dots, \vd{q}_{N_\omega - 1}^T]^T$ where $\vd{q}_j$ is the state at time $t_j = j \Delta T$. Throughout the paper, a subscript $\set{J}$ is used to indicate the set or vector that is formed by taking the subscript to be each element of the vector $\set{J} = [0, 1, \dots, N_\omega - 1]^T$. For example $t_\set{J}$ denotes the vector $[0,\Delta t, \dots, (N_\omega - 1)\Delta t]^T$. We denote the number of time steps $N_\omega$ for reasons related to the discrete Fourier transform, which will become clear later.
\\

An encoding scheme for trajectories comprises an encoder $\vd{e}^{x,t}: \mathbb{C}^{N_xN_\omega} \to \mathbb{C}^{rN_\omega}$ and a decoder $\vd{d}^{x,t}_\set{J}: \mathbb{C}^{rN_\omega} \to \mathbb{C}^{N_xN_\omega}$. The encoding scheme is effective if $\vd{q}_\set{J} \approx \vd{d}^{x,t}_\set{J}(\vd{e}^{x,t}(\vd{q}_\set{J}))$ for trajectories of the system. More precisely, the expected error averaged over the trajectory, $\mathbb{E}[ \| \vd{q}_\set{J} - \vd{d}^{x,t}_\set{J}(\vd{e}^{x,t}(\vd{q}_\set{J})) \|_{x,t}^2]$, should be small. The space-time norm is induced in the usual way from a space-time inner product $\| \vd{q}_\set{J}\|_{x,t} = \sqrt{ \langle \vd{q}_\set{J}, \vd{q}_\set{J} \rangle_{x,t} }$, and this inner product is given by
\begin{equation} \label{eq:trajrep:stip}
    \langle \vd{q}^1_\set{J}, \vd{q}^2_\set{J} \rangle_{x,t}  = \sum_{j = 0}^{N_\omega - 1} \langle \vd{q}^1_j , \vd{q}^2_j \rangle_{x} \text{.}
\end{equation}
\\

The analog of POD for the case of trajectories is space-time POD \cite{Lumley67,Frame23} --- it is the optimal linear encoding / decoding scheme for trajectories in the sense that it minimizes $\mathbb{E}[ \| \vd{q}_\set{J} - \vd{d}^{x,t}_\set{J}(\vd{e}^{x,t}(\vd{q}_\set{J})) \|_{x,t}^2]$. As discussed in the following subsection, SPOD modes can be used to represent trajectories, and a primary motivation for this paper is that the SPOD representation approaches the accuracy of the space-time POD representation for statistically stationary systems as $T \to \infty$. Another important point of comparison is the POD representation of trajectories, wherein each time step of the trajectory is encoded using POD modes, i.e., where $\vd{d}^{x,t}_j(\vd{e}^{x,t}(\vd{q}_\set{J})) = \td{\Phi}\td{\Phi}^*\td{W} \vd{q}_j$.

\subsection{Spectral proper orthogonal decomposition}
Spectral POD can be viewed as a `POD of the frequency domain' for statistically stationary flows: for every temporal frequency, there is a basis of SPOD modes that optimally represent the spatial structures exhibited at that frequency. To make this precise, we define the trajectory at frequency $\omega_k$ using the discrete Fourier transform (DFT) as
\begin{equation} \label{eq:SPOD:DFT}
    \hat{\vd{q}}_k = \text{DFT}_k\left[  \vd{q}_\set{J} \right] = \sum_{j = 0}^{N_\omega - 1} \vd{q}_j e^{-i\omega_k t_j} \text{,}
\end{equation}
where $\omega_k = 2\pi(k - N_\omega \Theta(k - N_\omega / 2) ) / T$, where $\Theta$ is $0$ if its argument is negative, and is $1$ otherwise. If the DFT is applied to a new trajectory of the same system (or a different interval of the same long trajectory), the Fourier components will likely be different --- if all $\hat{\vd{q}}_k$ were the same, the trajectories would be identical. The SPOD modes optimally capture this variability of $\hat{\vd{q}}_k$ just as the POD modes capture the variability of $\vd{q}$. Specifically, the SPOD modes at frequency $\omega_k$ can be defined to maximize the statistical energy captured at frequency $\omega_k$, which is defined as
\begin{equation} \label{eq:SPOD:opimization}
     \lambda_k(\vd{\psi}) =  \frac{\mathbb{E}\big[ | \langle \vdh{q}_k, \vd{\psi} \rangle_{x} |^2 \big]}{\| \vd{\psi}\|_{x}^2 }  \text{.}
\end{equation}
The first SPOD mode at frequency $\omega_k$ is defined to maximize \eqref{eq:SPOD:opimization}, and the latter modes are defined to maximize the same function subject to orthogonality with respect to previous modes,
\begin{equation}
\begin{aligned}
        &\vd{\psi}_{km} = \argmax \lambda_k(\vd{\psi}) \quad \\ \text{subject to} \quad &\langle \vd{\psi}_{km} , \vd{\psi}_{kn} \rangle_{x} = \delta_{mn} \quad \text{for} \quad 1 \leq m \leq n \text{,}
\end{aligned}
\end{equation}
where $\vd{\psi}_{km}$ denotes the $m$-th mode at the $k$-th frequency. We denote the energies of the modes as $\lambda_{km} = \lambda_k(\vd{\psi}_{km})$. Note the similarity to (space-only) proper orthogonal decomposition: the SPOD modes at $\omega_k$ are to $\vdh{q}_k$ as the POD modes are to $\vd{q}$. 
\\

With the modes defined, $\hat{\vd{q}}_k$ may be represented approximately by a linear combination of the first $r_k$ modes,
\begin{equation} 
    \hat{\vd{q}}_k \approx \td{\Psi}_k \vd{a}_k \text{.}
\end{equation}
Here, $\td{\Psi}_k = [\vd{\psi}_{k1}, \vd{\psi}_{k2}, \dots, \vd{\psi}_{kr_k}] \in \mathbb{C}^{N_x \times r_k}$ is the truncated matrix of modes at frequency $\omega_k$ and $\vd{a}_k \in \mathbb{C}^{r_k}$ is the vector of expansion coefficients. It will become clear below why it is advantageous to retain different numbers of modes at different frequencies. The error in the representation can be expressed in terms of the truncated SPOD energies as
\begin{equation} \label{eq:SPOD:ex_err_w}
    \mathbb{E}[ \| \vdh{q}_k - \SPOD{k} \vd{a}_k \|_{x}^2] = \sum_{m = r_k + 1}^{N_x} \lambda_{km} \text{,}
\end{equation}
and the first $r_k$ SPOD modes minimize this quantity over all rank-$r_k$ bases. If $\hat{\vd{q}}_k$ is known, the SPOD coefficients may be found using
\begin{equation} \label{eq:SPOD:coefs_from_qhat}
    \vd{a}_k = \td{\Psi}_k^{*} \td{W} \hat{\vd{q}}_k \text{.}
\end{equation}
This may be written in terms of a space-time encoder $\vd{a}_k = \vd{e}^{x,t}_k(\vd{q}_\set{J})$, where the $k$-th frequency of the encoder is given by
\begin{equation}
    \vd{e}^{x,t}_k(\vd{q}_\set{J}) = \td{\Psi}_k^{*} \td{W} \ \text{DFT}_k[\vd{q}_\set{J}] \text{.}
\end{equation}
\\

To obtain the SPOD modes at frequency $\omega_k$ from data, a set of realizations of $\vdh{q}_k$ is needed. These realizations could come from separate runs of the system, or, more commonly, from possibly overlapping blocks of a single long run \cite{Welch67}. Grouping $N_d$ of these realizations into a data matrix $\tdh{Q}_k = [\vdh{q}_k^1, \vdh{q}_k^2,\dots, \vdh{q}_k^{N_d}] \in \mathbb{C}^{N_x \times N_d}$, the modes can be obtained by first taking a singular value decomposition (SVD) of the weighted data matrix $\frac{1}{N_d}\td{W}^{1/2} \tdh{Q}_k = \td{U}_k \td{\Sigma}_k \td{V}_k^*$. The first $N_d$ SPOD modes are then given by $\td{\Psi}_k^{N_d} = \td{W}^{-1/2}\td{U}_k \in \mathbb{C}^{N_x \times N_d}$. Likewise, the diagonal matrix of the first $N_d$ SPOD energies is given by the squares of the singular values, $\td{\Lambda}_k^{N_d} = \td{\Sigma}_k^2$.

\subsection{SPOD modes for trajectory representation}
Our interest in using SPOD modes stems from the fact that they provide a near-optimal representation of trajectories. The SPOD modes at each frequency form an approximate representation of the trajectory in the frequency domain. Accordingly, they can be used to give an approximate space-time representation of the trajectory in the time domain by taking an inverse DFT,
\begin{equation} \label{eq:SPOD:SPOD_st_rep}
    \tilde{\vd{q}}_j= \frac{1}{N_\omega}\sum_{k = 0}^{N_\omega-1} \SPOD{k} \vd{a}_k e^{i\omega_k t_j} \text{,}
\end{equation}
where $\vdt{q}_\set{J}$ is the approximate trajectory. It may be written in terms of a space-time decoder as $\vdt{q}_\set{J} = \vd{d}^{x,t}_\set{J}(\vd{a}_\set{K})$, where the $j$-th time of the decoder is given by
\begin{equation}
    \vd{d}^{x,t}_{j}(\vd{a}_\set{K}) = \frac{1}{N_\omega}\sum_{k = 0}^{N_\omega-1} \SPOD{k} \vd{a}_k e^{i\omega_k t_j} \text{.}
\end{equation}
Hereafter, we use a subscript $\set{K}$ to denote all frequencies, analogous to our use of the subscript $\set{J}$ for all times. The relevant metric for accuracy is the spatial error averaged over the times in the interval, which is the square of the space-time norm of the difference. Due to Parseval's theorem, this is proportional to the error summed over frequencies,
\begin{equation} \label{eq:SPOD:sterror}
    \| \vdt{q}_\set{J} - \vd{q}_\set{J} \|^2_{x,t} = \frac{1}{N_\omega}\sum_{k = 0}^{N_\omega -1} \| \SPOD{k}\vd{a}_k - \vdh{q}_k \|_{x}^2 \text{.}
\end{equation}
Taking the expected value of both sides and using \eqref{eq:SPOD:ex_err_w}, we have
\begin{equation}
    \mathbb{E} \left[  \| \vdt{q}_\set{J} - \vd{q}_\set{J} \|^2_{x,t} \right] = \frac{1}{N_\omega} \sum_{k = 0}^{N_\omega -1} \sum_{m = r_k + 1}^{N_x} \lambda_{km} \text{.}
\end{equation}
Thus, the expected accuracy of the trajectory representation is determined by the energies of the truncated SPOD modes. This informs the best allocation of modes across frequencies: for the best accuracy given a total number of modes, one should retain the modes with the highest energies over all frequencies. We denote the average number of modes retained by $r$, and denote the $j$-th largest SPOD energy at any frequency by $\tilde{\lambda}_j$. The number of modes retained at the $k$-th frequency $r_k$ is then the number of modes at this frequency with energy greater than or equal to $\tilde{\lambda}_{(rN_\omega)}$,
\begin{equation} \label{eq:SPOD:rk}
    r_k = |\{ m: \lambda_{km} \geq \tilde{\lambda}_{(rN_\omega)} \}| \text{,}
\end{equation}
where $| \cdot |$ denotes the cardinality of the set, and the parentheses in the subscript are a reminder that $r$ and $N_\omega$ are multiplied and are not separate indices. 
\\

The SPOD trajectory representation \eqref{eq:SPOD:SPOD_st_rep} may be viewed as a rank-$rN_\omega$ space-time representation, i.e., it is a sum of spatiotemporal modes with constant coefficients of the form $\sum_{i = 1}^{rN_\omega} \vd{\xi}_i a_i$ with $\vd{\xi}_i \in \mathbb{C}^{N_\omega N_x}$. In the SPOD representation, each spatiotemporal mode is an SPOD mode at some frequency $\omega_k$ multiplied by the associated oscillatory time dependence $e^{i\omega_k t}$. The most accurate space-time representation is space-time POD modes --- these lead to lower trajectory representation error than any other space-time basis. The fact that motivates our use of SPOD modes is that as the time interval $T \to \infty$, the SPOD modes approach space-time POD modes \cite{Lumley67,Lumley70,Frame23}. This limit is approached fairly quickly in practice, so when $T$ is long compared to the correlation time of the system, the SPOD modes provide a representation of trajectories that is nearly optimally accurate given the total number of coefficients used (see Figure 2 in Ref. \cite{Frame25}). 
\\

\begin{figure}[]
    \centering
    \includegraphics{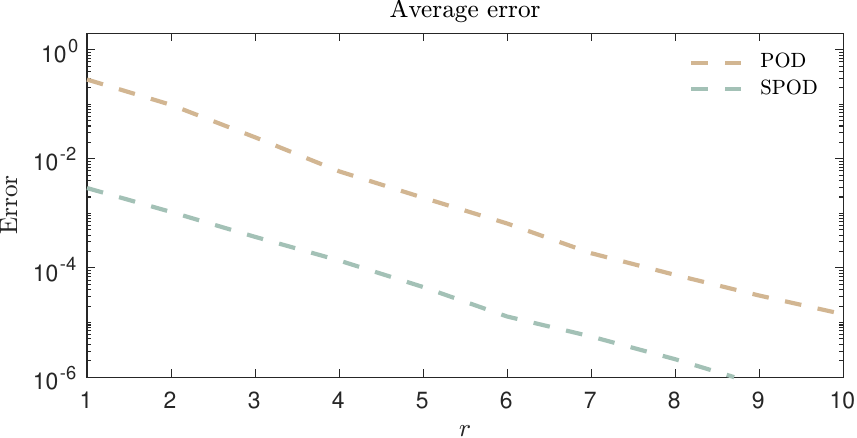}
    \caption{The relative error for the SPOD and POD trajectory encoders as a function of $r$, the number of modes used.}
    \label{fig:rep_err}
\end{figure}
The efficacy of the SPOD modes in representing trajectories can be compared to that of POD modes by calculating the relative error produced by both encoding/decoding. The relative error is 
\begin{equation}
    \frac{\|\vdt{q}_\set{J} - \vd{q}_\set{J}\|_{x,t}^2} {\|\vd{q}_\set{J}\|^2_{x,t}} \text{,}
\end{equation}
where $\vdt{q}_\set{J}$ is given by either the POD or SPOD encoding/decoding. In Figure~\ref{fig:rep_err}, we show these errors as a function of $r$, the number of modes, averaged over $30$ trajectories for the cubically nonlinear Ginzburg-Landau system described later in the paper. Figure~\ref{fig:rep_err} demonstrates the superior representational power of SPOD modes over POD modes; with the same number of total coefficients, the SPOD representation results in more than two orders of magnitude lower error. The goal of the paper is to develop a reduced-order model for nonlinear systems that solves for the SPOD coefficients quickly and accurately, so as to recover the accuracy shown in Figure~\ref{fig:rep_err} at low cost.

\subsection{Linear spectral solution operator projection} \label{sec:Lin}
Here, we review our previous work on SSOP in the linear setting. The full details are given in Ref. \cite{Frame25}. 
\\

The method aims to find, as accurately as possible, the SPOD coefficients of the solution $\vd{q}_\set{J}$ (defined on a temporal grid) to the LTI system
\begin{equation} \label{eq:Lin:ODE}
    \dot{\vd{q}} = \td{A}(\vd{\mu})\vd{q} + \td{B}\vd{f} \text{,}
\end{equation}
where $\vd{f} \in \mathbb{C}^{N_f}$ is a known time-dependent forcing. $\td{A}(\vd{\mu}) \in \mathbb{C}^{N_x \times N_x}$ and $\td{B} \in \mathbb{C}^{N_x \times N_f}$ are the usual system matrices, the former of which depends affinely on the parameter vector $\vd{\mu} \in \mathbb{R}^{N_\mu}$. We suppress this parameter dependence for now; after deriving the nonlinear reduced-order model in Section~\ref{sec:method}, we reintroduce it and show how the reduced matrices involved may be computed quickly given $\vd{\mu}$.
\\

The method derives from the fact that the SPOD coefficients are given by linear functionals of the initial condition and forcing. Specifically, the $m$-th SPOD coefficient at the $k$-th frequency is 
\begin{equation} \label{eq:Lin:ipSSOP}
    a_{km} = \langle \vd{L}[\vd{q}_0, \td{B}\vd{f}(t); t_\set{J}], \vd{\psi}_{km} e^{i\omega_k t_\set{J}} \rangle_{x,t} \text{,}
\end{equation}
where $\vd{a}_k=[a_{k1}, \dots , a_{kr_k}]^T$, and where the space-time inner product is \eqref{eq:trajrep:stip}. $\vd{L}[\vd{q}_0, \vd{g}(t); t_j]$ is a solution operator that maps the initial condition, forcing, and a time $t_j$ to $\vd{q}(t_j)$ as
\begin{equation} \label{eq:Lin:solop}
    \vd{L}[\vd{q}_0, \vd{g}(t); t_j] = e^{\td{A}t_j}\vd{q}_0 + \int_0^{t_j} e^{\td{A}(t_j - t')} \vd{g}(t') \ \text{d}t' \text{.} 
\end{equation}
\\

Analytic progress can be made by breaking \eqref{eq:Lin:ipSSOP} into a DFT of the solution operator and a spatial inner product as
\begin{subequations}
\begin{equation} \label{eq:Lin:DFT_plugin}
    \vdh{q}_k = \sum_{j = 0}^{N_\omega - 1} \left( e^{\td{A}j\Delta t}\vd{q}_0 + \int_0^{j \Delta t} e^{\td{A} (j \Delta t - t')}  \td{B} \vd{f}(t') \ dt' \right)  e^{-i\omega_k j \Delta t} \text{.}
\end{equation}
\begin{equation} \label{}
  \vd{a}_k = \td{\Psi}_{k}^* \td{W} \vdh{q}_k
\end{equation}
\end{subequations}
where we the left multiplication by $ \td{\Psi}_{k}^* \td{W}$ effects the inner product for all modes at the $k$-th frequency. The sum in \eqref{eq:Lin:DFT_plugin} can be computed analytically assuming that \textit{(i)} the forcing is of the form $\vd{f}(t) = \sum_{l = 0}^{N_\omega-1} \vdh{f}_l e^{i\omega_lt}$, and \textit{(ii)}  that $i\omega_l \td{I} - \td{A}$ is invertible for all included $l$ (this latter assumption can be relaxed). Under these assumptions, \eqref{eq:Lin:DFT_plugin} can be written
\begin{equation} \label{eq:Lin:qhat(f,q0)}
    \hat{\vd{q}}_k = \td{R}_k \td{B}  \hat{\vd{f}}_k + \left(\td{I} -e^{(\td{A} - i \omega_k \td{I})\Delta t} \right)^{-1} \left(\td{I} - e^{\td{A}T} \right) \left(\vd{q}_0 - \frac{1}{N_\omega} \sum_{l = 0}^{N_\omega - 1} \td{R}_l \td{B} \hat{\vd{f}}_l \right) \text{.}
\end{equation}
Here, $\td{R}_k = (i\omega \td{I} - \td{A})^{-1}$ is the resolvent operator. Though assumption \textit{(i)} will seldom hold in practice, we have found in our numerical experiments that this introduces minimal error. Using \eqref{eq:SPOD:coefs_from_qhat}, we obtain an equation for the SPOD coefficients,
\begin{equation} \label{eq:Lin:a(f)_exact}
     \vd{a}_k = \SPOD{k}^{*} \td{W}\td{R}_k \td{B}  \hat{\vd{f}}_k + \SPOD{k}^{*} \td{W}\left(\td{I} -e^{(\td{A} - i \omega_k \td{I})\Delta t} \right)^{-1} \left(\td{I} - e^{\td{A}T} \right) \left(\vd{q}_0 - \frac{1}{N_\omega} \sum_{l = 0}^{N_\omega - 1} \td{R}_l \td{B} \hat{\vd{f}}_l \right) \text{.}
\end{equation}
\\

At this point, a number of approximations must be made to make the operators in \eqref{eq:Lin:a(f)_exact} tractable to compute for large systems, i.e., to reduce the offline cost, and to make the equation fast to evaluate for all frequencies online. If the system in question is small enough for direct computation of the matrices in \eqref{eq:Lin:a(f)_exact}, then only the last approximation described in this subsection is necessary. 
\\

First, the operator $ \SPOD{k}^{*} \td{W}\td{R}_k \td{B}$ must be approximated. We do this by making use of the data matrices used to compute the SPOD modes to obtain the action of the resolvent operator within a subspace, then form an approximation by projecting into this subspace. Specifically, we form the matrix $\td{G}_k = \td{L}_k\tdh{Q}_k \in \mathbb{C}^{N_x \times N_d}$, where  $\td{L}_k = i\omega_k\td{I} - \td{A}$. This implies
\begin{equation} \label{eq:Lin:Gdef}
    \tdh{Q}_k = \td{R}_k \td{G}_k \text{,}   
\end{equation}
meaning we have the action of the resolvent operator on the subspace spanned by the columns of $\td{G}_k$. Using \eqref{eq:Lin:Gdef}, we approximate the resolvent using
\begin{equation} \label{eq:Lin:Raprx}
    \td{R}_k \approx \tdh{Q}_k\td{G}_k^+ \text{,}
\end{equation}
where $\td{G}_k^+ = (\td{G}_k^* \td{W} \td{G}_k)^{-1} \td{G}_k^* \td{W}$ is the Moore-Penrose pseudoinverse of $\td{G}_k$. It may be shown of this approximation that
\begin{equation} \label{eq:Lin:Raprx_proj}
   \tdh{Q}_k\td{G}_k^+ = \td{R}_k \td{P}_{G_k}\text{,}
\end{equation}
where $\td{P}_{G_k}= \td{G}_k(\td{G}_k^* \td{W} \td{G}_k)^{-1} \td{G}_k^* \td{W}$ is the orthogonal projection operator into the column space of $\td{G}_k$. This implies that the approximation \eqref{eq:Lin:Raprx} of the resolvent is accurate to the extent that it is applied to vectors near the column space of $\td{G}_k$. The approximation of $ \SPOD{k}^{*} \td{W}\td{R}_k \td{B}$ is then 
\begin{equation}\label{eq:Lin:Edef}
      \td{E}_k = \SPOD{k}^* \td{W} \tdh{Q}_k \td{G}^+ \td{B}  \approx \SPOD{k}^{*} \td{W}\td{R}_k \td{B} \text{.}
\end{equation}
\\

Next, we approximate $\SPOD{k}^{*} \td{W}\left(\td{I} -e^{(\td{A} - i \omega_k \td{I})\Delta t} \right)^{-1} \left(\td{I} - e^{\td{A}T} \right)$. This is accomplished by the Galerkin approximation
\begin{equation}
    \td{P}_k \left(\td{I} -e^{(\tilde{\td{A}} - i \omega_k \td{I})\Delta t} \right)^{-1} \left(\td{I} - e^{\tilde{\td{A}}T} \right) \SPOD{k}^{N_d*} \td{W} \approx \SPOD{k}^{*} \td{W}\left(\td{I} -e^{(\td{A} - i \omega_k \td{I})\Delta t} \right)^{-1} \left(\td{I} - e^{\td{A}T} \right)\text{,}
\end{equation}
where $\tilde{\td{A}} = \SPOD{k}^{N_d *}\td{W}\td{A} \SPOD{k}^{N_d} \in \mathbb{C}^{N_d \times N_d}$ and $\td{P}_k = \begin{bmatrix} \td{I} & \td{0} \end{bmatrix} \in \mathbb{R}^{r_k \times N_d}$ selects the first $r_k$ rows of the matrix it operates on. In Appendix~\ref{app:approxH}, we detail a more accurate, but more costly, approximation of $\SPOD{k}^{*} \td{W}\left(\td{I} -e^{(\td{A} - i \omega_k \td{I})\Delta t} \right)^{-1} \left(\td{I} - e^{\td{A}T} \right)$.
\\

The purpose of the previous two approximations is to alleviate the offline cost of building the operators in \eqref{eq:Lin:a(f)_exact}. In order to avoid poor online scaling, we approximate the term $ \vd{q}_0 - \frac{1}{N_\omega} \sum_{l = 0}^{N_\omega - 1} \td{R}_l \td{B} \hat{\vd{f}}_l $ by representing it in a rank-$p_1$ reduced spatial basis $\td{\Phi} \in \mathbb{C}^{N_x \times p}$. We refer to $\td{\Phi}$ as the intermediary basis, and we take it to be the (space-only) POD modes. This is helpful because with it, evaluating the impact of $ \vd{q}_0 - \frac{1}{N_\omega} \sum_{l = 0}^{N_\omega - 1} \td{R}_l \td{B} \hat{\vd{f}}_l $ for each frequency scales with $p_1$ rather than with $N_x$. The approximation is 
\begin{equation}
      \td{\Phi} \td{\Phi}^* \vd{q}_0 - \frac{1}{N_\omega} \sum_{l = 0}^{N_\omega - 1} \td{\Phi} \td{J}_l \hat{\vd{f}}_l  \approx \vd{q}_0 - \frac{1}{N_\omega} \sum_{l = 0}^{N_\omega - 1} \td{R}_l \td{B} \hat{\vd{f}}_l \text{,}
\end{equation}
where $\td{J}_k \in \mathbb{C}^{p_1 \times N_f}$ is an approximation of $\td{\Phi}^* \td{W} \td{R}_k \td{B}$. Using the same approximation of the resolvent from before, it is defined as
\begin{equation} \label{eq:Lin:Jdef}
    \td{J}_k = \td{\Phi}^* \td{W} \tdh{Q}_k \td{G}^+_k \td{B} \text{.}
\end{equation}
We also define the matrix $\td{H}_k \in \mathbb{C}^{r_k \times p_1}$ as
\begin{equation} \label{eq:Lin:Hdef} 
 \td{H}_k =  \td{P}_k \left(\td{I} -e^{(\tilde{\td{A}} - i \omega_k \td{I})\Delta t} \right)^{-1} \left(\td{I} - e^{\tilde{\td{A}}T} \right) \SPOD{k}^{N_d*} \td{W} \td{\Phi}  \text{.}
\end{equation}
The ROM equations are then 
\begin{equation} \label{eq:Lin:ROMfinal}
     \vdt{a}_k = \td{E}_k \vdh{f}_k + \td{H}_k \left( \td{\Phi}^* \vd{q}_0 - \frac{1}{N_\omega} \sum_{l = 0}^{N_\omega - 1}  \td{J}_l \hat{\vd{f}}_l \right) \text{,}
\end{equation}
where the $\tilde{(\cdot)}$ indicates the coefficients above are an approximation of the SPOD coefficients for the trajectory --- they are approximate due to the operator approximations introduced above. In \eqref{eq:Lin:ROMfinal}, the term in parentheses need only be computed once, and the online calculations scale as $\mathcal{O}\big(N_\omega (rN_f + rp_1 + N_f \log N_\omega)\big)$. Pseudocode for the offline and online phases of the method can be found in Ref.~\cite{Frame25}. The pseudocode given in Appendix~\ref{app:algorithms} can also be used for the linear problem by disregarding the nonlinearity.
\\

\section{Nonlinear spectral solution operator projection} \label{sec:method}
We now consider systems of the form
\begin{equation} \label{eq:NL_method:dy_sys}
    \dot{\vd{q}} = \td{A}\vd{q} + \td{B} \vd{f} + \vd{n}(\vd{q}) \text{,}
\end{equation}
where $\vd{n}: \mathbb{C}^{N_x} \to \mathbb{C}^{N_x}$ is some nonlinear function of the state. The parameter dependence in the linear operator $\td{A}$ is temporarily suppressed, and we assume the nonlinear operator is not parameter-dependent. We again assume that the system has been transformed such that the solution has zero mean. An additional transformation may be useful if the linear operator has unstable eigenvalues. The transformation is $\td{A} \to \td{A} - \alpha \td{I}$, $\vd{n}(\vd{q}) \to \vd{n}(\vd{q}) + \alpha \vd{q}$, where $\alpha$ is to the right of the spectrum of $\td{A}$, i.e., $\alpha > \max \Re( \lambda(\td{A})) $. Without this transformation, the terms related to the linear operator may dominate the solution.
\\

Our aim is to generalize SSOP to nonlinear systems of the form \eqref{eq:NL_method:dy_sys}. This is accomplished by treating the nonlinear term as an additional forcing on the system, but one that depends on the trajectory. The linear ROM \eqref{eq:Lin:ROMfinal} is amended by adding a term corresponding to the nonlinearity at frequency $\omega_k$ that results from the SPOD coefficients at all frequencies. As we show in Section~\ref{sec:nonlinear:implicit}, this results in a system that must be solved for all the SPOD coefficients at once. We then discuss two methods for approximating the effect of nonlinearity in Section~\ref{sec:nonlinear:construct}. Affine parameter dependence is discussed in Section~\ref{sec:nonlinear:parameter}. We present a method for solving the nonlinear system of equations in Section~\ref{sec:nonlinear:solve}, and discuss scaling in Section~\ref{sec:nonlinear:scaling}. 
\\

\begin{figure}[ht!]
    \centering
    \includegraphics[scale=0.6]{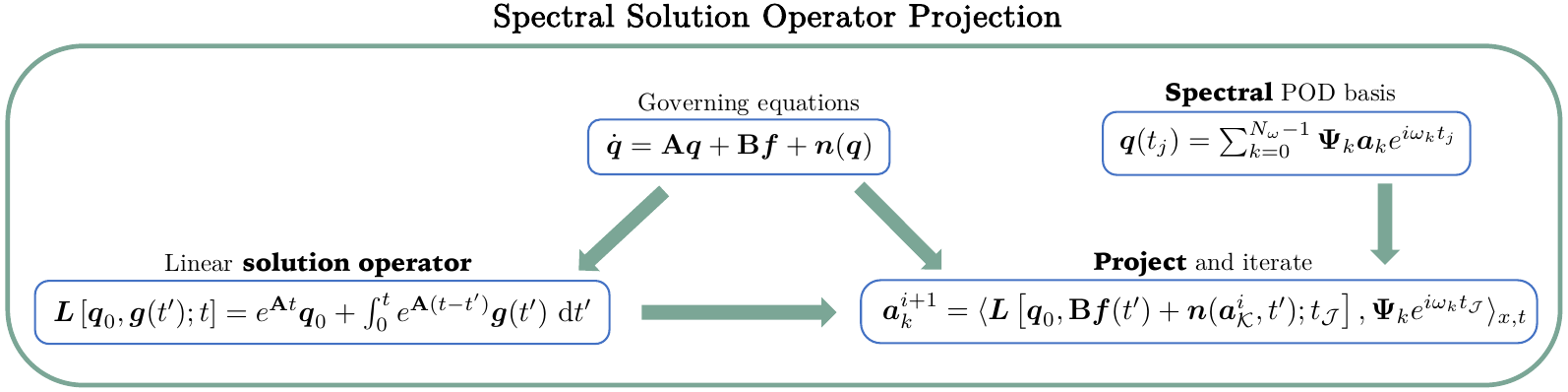}
    \caption{\textbf{Spectral Solution Operator Projection}: A system of nonlinear algebraic equations is obtained by performing a space-time projection of a linear solution operator onto the set of SPOD modes used to encode trajectories. The linear solution operator returns the solution to the linear system given the initial condition and forcing. Nonlinearities are handled in this framework as a state-dependent forcing. The system is solved via a fixed-point iteration, which may be viewed as a perturbation series around the linear solution.}
    \label{fig:Graphical_abstract}
\end{figure}

The method is shown schematically in Figure~\ref{fig:Graphical_abstract}, and the details of the components therein will be discussed in the remainder of Section~\ref{sec:method}. We note here that in the bottom-right box in the figure, an inner product is equated to a vector, and the second argument of the inner product is a matrix. In this shorthand notation, the $m$-th element of the vector is the inner product of the first argument with the $m$-th column of the matrix. 
\\

\subsection{Formulation} \label{sec:nonlinear:implicit}
To handle nonlinearity within the framework described above, we begin with the implicit equation
\begin{equation}
    \vd{q}_j = \vd{L}\left[ \vd{q}_0, \td{B}\vd{f}(t) + \vd{n}(\vd{q}(t)); t_j \right] \text{,}
\end{equation}
where $\vd{L}$ is the solution operator given in \eqref{eq:Lin:solop}. This formulation treats the nonlinearity as an additional forcing on the linear part of the system, but, crucially, one that depends on the state. The $m$-th SPOD coefficient at the $k$-th frequency is now
\begin{equation} \label{eq:nonlinear:SSOP}
    a_{km} = \langle \vd{L}\left[ \vd{q}_0, \td{B}\vd{f}(t) + \vd{n}(\vd{q}(t)); t_\set{J} \right] , \vd{\psi}_{km}e^{i\omega_k t_\set{J}} \rangle_{x,t} \text{.}
\end{equation}
In order to arrive at a closed system of equations, $\vd{n}(\vd{q}(t))$ must be approximated using the SPOD coefficients, and we return to this closure in Section~\ref{sec:nonlinear:construct}. First, as in the linear case, we evaluate the space-time inner product by taking a temporal DFT of the solution operator, followed by a spatial inner product:
\begin{subequations}
\begin{equation} \label{eq:nonlinear:DFT}
    \vdh{q}_k = \sum_{j = 0}^{N_\omega - 1} \left( e^{\td{A}j\Delta t}\vd{q}_0 + \int_0^{t_j} e^{\td{A} (t_j - t')}  \td{B} \vd{f}(t') + \vd{n}(\vd{q}(t')) \ dt' \right)  e^{-i\omega_k t_j} \text{.}
\end{equation}
\begin{equation} \label{eq:nonlinear:spatial_ip}
    \vd{a}_k = \td{\Psi}_k^* \td{W} \vdh{q}_k \text{.}
\end{equation}
\end{subequations}
With the two assumptions from Section~\ref{sec:Lin}, as well as the additional assumption \textit{(iii)} that the nonlinearity can be written in the form $\vd{n}(\vd{q} (t)) = \sum_{l = 0}^{N_\omega-1} \vdh{n}_l e^{i\omega_lt}$ where $\vdh{n}_l \in \mathbb{C}^{N_x}$, we show in Appendix~\ref{app:full_order_fd} that the $k$-th frequency of the solution is given by
\begin{equation} \label{eq:NL_method:qhat(f,q0)}
\begin{split} 
    \hat{\vd{q}}_k = &\td{R}_k \left( \td{B}  \hat{\vd{f}}_k + \vdh{n}_k^q(\vdh{q}_\set{K}) \right) \\ + &\left(\td{I} -e^{(\td{A} - i \omega_k\td{I})\Delta t} \right)^{-1} \left(\td{I} - e^{\td{A}T} \right) \left(\vd{q}_0 - \frac{1}{N_\omega} \sum_{l = 0}^{N_\omega - 1} \td{R}_l \left( \td{B}  \hat{\vd{f}}_l + \vdh{n}^q_l(\vd{q}_\set{K}) \right) \right) \text{.}
\end{split} 
\end{equation}
Here, $\vdh{n}^q_k: \mathbb{C}^{N_xN_\omega} \to \mathbb{C}^{N_x}$ returns the nonlinearity at frequency $\omega_k$ that results from $\vdh{q}_\set{K}$, the concatenation of all frequencies of the trajectory. The superscript is used to distinguish this function from a similar one introduced in the next subsection. One way to compute $\vdh{n}^q_k$ is to first take the inverse DFT of $\vdh{q}_\set{K}$, giving the trajectory in the time domain $\vd{q}_\set{J}$. The $k$-th frequency of the nonlinearity $\vdh{n}^q_k$ is given by then computing $\vd{n}(\vd{q}_\set{J})$ and taking the $k$-th component of the DFT of the result. Similar to the assumption of the form of the forcing, assumption \textit{(iii)} will likely not hold in practice, but we have found that the error introduced is minimal.
\\

Inserting \eqref{eq:NL_method:qhat(f,q0)} into \eqref{eq:nonlinear:spatial_ip}, we have
\begin{equation} \label{eq:nonlinear:first_spatial_red}
\begin{split} 
    \vdt{a}_k = &\SPOD{k}^{*} \td{W} \td{R}_k \left( \td{B}  \hat{\vd{f}}_k + \vdh{n}_k(\vdt{a}_\set{K}) \right) \\ + &\SPOD{k}^{*} \td{W}  \left(\td{I} -e^{(\td{A} - i \omega_k)\Delta t} \right)^{-1} \left(\td{I} - e^{\td{A}T} \right) \left(\vd{q}_0 - \frac{1}{N_\omega} \sum_{l = 0}^{N_\omega - 1} \td{R}_l \left( \td{B}  \hat{\vd{f}}_l + \vdh{n}_l(\vdt{a}_\set{K}) \right) \right) \text{.}
\end{split}
\end{equation}
Here, $\vdt{a}_\set{K}$ is the (ROM approximation of the) vector of retained SPOD coefficients at all frequencies concatenated together. With the reduction, the nonlinearity at each frequency must be computed as a function of the vector $\vdt{a}_\set{K}$, not $\vdh{q}$, i.e., $\vdh{n}_k: \mathbb{C}^{rN\omega} \to \mathbb{C}^{N_x}$. Two different methods for approximating the nonlinearity at each frequency, given the SPOD coefficients, are described in the following subsection. 

The operators in \eqref{eq:nonlinear:first_spatial_red} must be approximated if the system is large. With the same approximations used in the linear case, we have,
\begin{equation} \label{eq:nonlinear:spatial_red_lin_apprxs}
     \vdt{a}_k = \td{E}_k \vdh{f}_k + \SPOD{k}^{*} \td{W} \td{R}_k  \vdh{n}_k(\vdt{a}_\set{K}) + \td{H}_k \left( \td{\Phi}^* \vd{q}_0 - \sum_{l = 0}^{N_\omega - 1} \td{J}_l \vdh{f}_l + \td{\Phi}^* \td{W} \td{R}_k  \vdh{n}_k(\vdt{a}_\set{K}) \right) \text{.}
\end{equation}
\\

\subsection{Constructing the nonlinearity} \label{sec:nonlinear:construct}
To close the reduced-order model equations \eqref{eq:nonlinear:spatial_red_lin_apprxs}, we must compute the $k$-th frequency of the nonlinearity $\vdh{n}_k$ that results from the SPOD coefficients. With no approximation, the nonlinearity is given in terms of the $k$-th frequency of the nonlinear function acting on the decoding of the coefficients 
\begin{equation} \label{eq:nonlinear:n_k}
    \vdh{n}_k(\vdt{a}_\set{K}) = \text{DFT}_k \left[ \vd{n}(\vd{d}_\set{J}^{x,t} (\vdt{a}_\set{K} ) ) \right] = \text{DFT}_k \left[ \vd{n}\left( \sum_{l = 0}^{N_\omega - 1} \SPOD{l}\vdt{a}_l e^{i\omega_l t_\set{J}} \right) \right] \text{.}
\end{equation}
However, computing the nonlinearity this way is expensive since it involves expanding to the $N_x$-dimensional space. The problem of obtaining $\vdh{n}_k$ without doing computations in this space is entirely analogous to approximating the nonlinearity at a given time step given the POD coefficients at that time step. Accordingly, we use hyper-reduction, specifically, the discrete empirical interpolation method \cite{Chaturantabut10} (DEIM) in the case that the nonlinearity is not a quadratic form. If the nonlinearity is a quadratic form, we use (sparsified) triadic interactions. We note that it may be appropriate to apply a dealiasing technique \cite{Orszag71} in computing the nonlinearity, such as truncating the sum in \eqref{eq:nonlinear:n_k}. The details and effectiveness of dealiasing depend on the nature of the nonlinearity. We therefore ignore this issue in dealing with general nonlinearities, but implicitly employ a standard 2/3-rule dealiasing for the quadratic nonlinearity. We note here that the function $\vd{n}(\vdt{a}_\set{K},t)$ in Figure~\ref{fig:Graphical_abstract} is $\vdh{n}_\set{K}(\vdt{a}_\set{K})$ transformed to the time domain. 
\\

\subsubsection{Hyper-reduction} \label{sec:nonlinear:deim}
Here, we describe a DEIM-based approximation of $\vdh{n}_k(\vdt{a}_\set{K})$. We assume that the nonlinearity is a local function of the state, i.e., each component of the nonlinearity is a function of only the same component of the state $n_i(\vd{q}) = n(q_i)$. The approach can trivially be extended to cases where the nonlinearity is a local function of the state and some number of linear operations on it, such as $n_i(\vd{q}) = n(q_i,(\td{D}^1\vd{q})_i)$, where $\td{D}^1$ is a differentiation matrix.
\\

We build the matrices $\td{U}^n \in \mathbb{C}^{N_x \times p_2}$ and $\td{P}^{nT} \in \mathbb{R}^{p_2 \times N_x}$ following the standard DEIM algorithm using snapshots of the nonlinearity and $p_2$ sample points. The matrix $\td{U}^n$ is the POD modes for the nonlinearity and the matrix $\td{P}^{n} = [\vd{e}_{s_1}, \dots , \vd{e}_{s_{p_2}}]^T$ is a the sampling matrix, where $\vd{e}_j$ is the $j$-th canonical unit vector and $s_i$ is the $i$-th sample point given by DEIM. The standard DEIM approximation of the nonlinearity at time $t$ given the sampled state is $\vd{n}(\vd{q}) \approx \td{U}^n (\td{P}^{nT} \td{U}^n )^{-1} \vd{n}(\td{P}^{nT} \vd{q})$, where $\vd{n}(\td{P}^{nT} \vd{q}) \in \mathbb{C}^{p_2}$. 
\\

We adapt this methodology to approximating $\vdh{n}_k$ as follows. The nonlinearity at each time $t_j$ is approximated using DEIM and the state at $t_j$ that is implied by the SPOD coefficients. Then, $\vdh{n}_k$ is obtained using the DFT of the nonlinearity at all times as
\begin{equation} \label{eq:nonlinear:deim:aprx1} 
    \vdh{n}_k(\vdt{a}_\set{K}) \approx \text{DFT}_k \left[ \td{U}^n (\td{P}^{nT} \td{U}^n )^{-1} \td{P}^{nT} \vd{n}\left( \sum_{l = 0}^{N_\omega - 1} \SPOD{l}\vdt{a}_l e^{i\omega_l t_\set{J}} \right) \right] \text{.}
\end{equation}
The right-hand side of \eqref{eq:nonlinear:deim:aprx1} is expensive to evaluate as written. However, leveraging the locality of the nonlinearity, the sampling operator $\td{P}^{nT}$ may be applied directly to the SPOD modes, and the matrix $\td{U}^n (\td{P}^{nT} \td{U}^n )^{-1}$ may be moved outside the DFT, as it is independent of time. Then, we have
\begin{equation} \label{eq:nonlinear:deim:aprx2}
\begin{split}
    &\text{DFT}_k \left[ \td{U}^n (\td{P}^{nT} \td{U}^n )^{-1} \td{P}^{nT} \vd{n}\left( \sum_{l = 0}^{N_\omega - 1} \SPOD{l}\vdt{a}_l e^{i\omega_l t_\set{J}} \right) \right] \\ = &\td{U}^n (\td{P}^{nT} \td{U}^n )^{-1} \text{DFT}_k \left[ \vd{n}\left( \sum_{l = 0}^{N_\omega - 1} \td{P}^{nT} \SPOD{l}\vdt{a}_l e^{i\omega_l t_\set{J}} \right) \right] \text{.}
\end{split}
\end{equation}
With the appropriate precomputed operators, evaluating the expression on the right does not scale with $N_x$, as desired. Three operators must be precomputed. First, the sampling matrix applied to the SPOD modes, $\td{S}_k = \td{P}^{nT} \SPOD{k}^{}$, must be precomputed for each frequency. Next, the matrix $\SPOD{k}^{*} \td{W} \td{R}_k \td{U}^n (\td{P}^{nT} \td{U}^n )^{-1}$ must be both precomputed and approximated for each frequency. Using the approximation $\td{R}_k \approx  \tdh{Q}_k\td{G}_k^+$ established in Section~\ref{sec:Lin}, we have
\begin{equation} \label{eq:nonlinear:deim:matrix1}
    \td{N}_k = \SPOD{k}^{*} \td{W} \tdh{Q}_k\td{G}_k^+ \td{U}^n (\td{P}^{nT} \td{U}^n )^{-1} \approx \SPOD{k}^{*} \td{W} \td{R}_k \td{U}^n (\td{P}^{nT} \td{U}^n )^{-1} \in \mathbb{C}^{r_k \times p_2} \text{.}
\end{equation}
Similarly, the matrix $\td{\Phi}^* \td{W} \td{R}_k \td{U}^n (\td{P}^{nT} \td{U}^n )^{-1}$ must be precomputed and approximated at each frequency. Using the same approximation of the action of the resolvent, we have
\begin{equation}\label{eq:nonlinear:deim:matrix2}
     \td{M}_k = \td{\Phi}^* \td{W} \tdh{Q}_k\td{G}_k^+ \td{U}^n (\td{P}^{nT} \td{U}^n )^{-1} \approx \td{\Phi}^* \td{W} \td{R}_k \td{U}^n (\td{P}^{nT} \td{U}^n )^{-1} \in \mathbb{C}^{p_1 \times p_2} \text{.}
\end{equation}
\\

 With these approximations in place, the nonlinear system with the DEIM approximation of the nonlinearity is 
 \begin{equation} \label{eq:nonlinear:deim:system}
        \vdt{a}_k = \td{E}_k \vdh{f}_k + \td{N}_k \vdh{n}_k^s(\vdt{a}_\set{K}) + \td{H}_k \left( \td{\Phi}^* \vd{q}_0 - \sum_{l = 0}^{N_\omega - 1} \td{J}_l \vdh{f}_l + \td{M}_l \vdh{n}_l^s(\vdt{a}_\set{K}) \right) \text{,}
 \end{equation}
where the sampled nonlinearity $\vdh{n}_k^s \in \mathbb{C}^{p_2}$ is calculated as $\vdh{n}_k^s = \text{DFT}_k \left[ \vd{n}\left( \sum_{l = 0}^{N_\omega - 1} \td{S}_l \vdt{a}_l e^{i\omega_l t_\set{J}} \right) \right]$. 

\subsubsection{Sparse triadic interactions} \label{sec:nonlinear:bil}
Many systems of interest have quadratic nonlinearities, i.e., ones of the form $\vd{n}(\vd{q}) = \vd{b}(\vd{q},\vd{q})$, where $\vd{b}: \mathbb{C}^{N_x} \times \mathbb{C}^{N_x} \to \mathbb{C}^{N_x}$ is a bilinear function. In these systems, $\vdh{n}_k(\vdt{a}_\set{K})$ can be computed exactly in the reduced space by summing over triadic interactions. These sums can be costly, so we discuss a strategy for discarding many of these triadic interactions. The resulting approximation of $\vdh{n}_k(\vdt{a}_\set{K})$ is likely to be more accurate in many applications than the DEIM-based approximation described above.
\\

Expressing \eqref{eq:nonlinear:n_k} in terms of $\vd{b}$ gives
\begin{equation}
    \vdh{n}_k(\vdt{a}_\set{K}) = \text{DFT}_k \left[ \vd{b}\left( \sum_{l = 0}^{N_\omega - 1} \SPOD{l}\vdt{a}_l e^{i\omega_l t_\set{J}}, \sum_{i = 0}^{N_\omega - 1} \SPOD{i}\vdt{a}_i e^{i\omega_i t_\set{J}} \right) \right] \text{.}
\end{equation}
After leveraging the bilinearity, we have
\begin{equation} \label{eq:nonlinearity:bilinear_intermediate2}
    \vdh{n}_k(\vdt{a}_\set{K}) = \text{DFT}_k \left[\sum_{l = 0}^{N_\omega - 1} \sum_{i = 0}^{N_\omega - 1} e^{i(\omega_l + \omega_i)t_\set{J}} \vd{b}\left(  \SPOD{l}\vdt{a}_l ,  \SPOD{i}\vdt{a}_i \right) \right] \text{.}
\end{equation}
$\text{DFT}_k[\cdot]$ only depends on pairs of frequencies where $l + i \equiv k \ (\text{mod} \ N_\omega)$, however some of these terms are generated by aliasing. After expanding the SPOD mode sums, making use of the bilinearity, and dropping the aliasing terms, the nonlinear term becomes
\begin{equation} \label{eq:nonlinear:bilinear_intermediate3}
    \vdh{n}_k(\vdt{a}_\set{K}) = \sum_{\omega_l + \omega_i = \omega_k} \sum_{m = 1}^{r_l} \sum_{n = 1}^{r_i} a_{lm} a_{in} \vd{b}\left(  \vd{\psi}_{lm} , \vd{\psi}_{in} \right)  \text{.}
\end{equation}
The terms involving $\vdh{n}_k(\vdt{a}_\set{K})$ in the spatially reduced equation for the SPOD coefficients \eqref{eq:nonlinear:spatial_red_lin_apprxs} are $\SPOD{k}^{*} \td{W} \td{R}_k  \vdh{n}_k(\vdt{a}_\set{K})$ and $\td{\Phi}^* \td{W} \td{R}_k  \vdh{n}_k(\vdt{a}_\set{K})$. To evaluate these terms quickly online, we precompute the vectors
\begin{subequations} \label{eq:nonlinear:all}
\begin{equation} \label{eq:nonlinear:nklmndef}
    \vd{n}_{klmn} = \SPOD{k}^{*} \td{W} \tdh{Q}_k \td{G}_k^+\vd{b}\left(  \vd{\psi}_{lm} , \vd{\psi}_{in} \right) \in \mathbb{C}^{r_k} \text{,}
\end{equation}
\begin{equation} \label{eq:nonlinear:mklmndef}
    \vd{m}_{klmn} = \td{\Phi}^* \td{W} \tdh{Q}_k \td{G}_k^+\vd{b}\left(  \vd{\psi}_{lm} , \vd{\psi}_{in} \right) \in \mathbb{C}^{p_1} \text{,}
\end{equation}
\end{subequations}
where $\tdh{Q}_k \td{G}_k^+$ is the approximation of $\td{R}_k$ used elsewhere in the paper. In (\ref{eq:nonlinear:all}), we omit an `$i$' subscript because the frequency $\omega_i$ is determined by $\omega_k$ and $\omega_l$.
\\

In many physical systems, the great majority of the triadic interactions are of negligible magnitude \cite{Schmidt20}. Accordingly, many of the terms in the sum in \eqref{eq:nonlinear:bilinear_intermediate3} may have little impact on $\vdh{n}_k(\vdt{a}_\set{K})$, due to either $|a_{lm}a_{in}|$ or $\| \vd{b}\left(  \vd{\psi}_{lm}, \vd{\psi}_{in} \right) \|_{x}$ being small. In order to accelerate the computation of $\vdh{n}_k$, it is useful to ignore these irrelevant terms in the online stage of the method. We set a threshold $\epsilon$ and only retain terms for which
\begin{equation} \label{eq:nonlinear:bilinear:inclusion_criterion}
    \| \vd{n}_{klmn} \|_{x}^2 \lambda_{lm}\lambda_{in} > \epsilon \text{.}
\end{equation}
The quantity $\| \vd{n}_{klmn} \|_{x}^2 \lambda_{lm}\lambda_{in}$ is an easy-to-evaluate proxy for $\mathbb{E}[\| \vd{n}_{klmn}a_{lm}a_{in} \|_{x}^2]$, and the two are equivalent in the case that $\mathbb{E}[|a_{lm}a_{in}|^2] = \mathbb{E}[|a_{lm}|^2]\mathbb{E}[|a_{in}|^2]$. We denote the set of tuples $(l,m,n)$ that meet the criterion \eqref{eq:nonlinear:bilinear:inclusion_criterion} by
\begin{equation} \label{eq:nonlinear:Ndef}
    \set{N}_k = \{(l,m,n) :  \| \vd{n}_{klmn} \|_{x}^2 \lambda_{lm}\lambda_{in} > \epsilon \} \text{.}
\end{equation}
The nonlinear system with the nonlinearity approximated using sparse triadic interactions is 
\begin{equation} \label{eq:nonlinear:bilinear:system}
    \vdt{a}_k = \td{E}_k \vdh{f}_k + \sum_{(l,m,n) \in \set{N}_k} \vd{n}_{klmn} a_{lm}a_{in} + \td{H}_k \left( \td{\Phi}^* \vd{q}_0 - \sum_{l = 0}^{N_\omega - 1} \td{J}_l \vdh{f}_l +  \sum_{(p,m,n) \in \set{N}_l} \vd{m}_{lpmn} a_{lm}a_{in} \right) \text{,}
\end{equation}
where the index `$i$' is implied by $l$ via $\omega_i = \omega_k - \omega_l$. 
\\

A practical matter bears mentioning: the nonlinear term in a quadratically nonlinear full-order model will likely not be written in quadratic form. That is, the function $\vd{n}(\vd{q})$ will be available, but not the function $\vd{b}(\vd{q}_1,\vd{q}_2)$. In fact, the latter is not uniquely defined. Defining $\vd{b}$ to also be symmetric, i.e., $\vd{b}(\vd{q}_1,\vd{q}_2) = \vd{b}(\vd{q}_2,\vd{q}_1) \ \forall \vd{q}_1,\vd{q}_2$, makes it both unique and easy to calculate from $\vd{n}$. The unique symmetric bilinear function with $\vd{b}(\vd{q},\vd{q}) = \vd{n}(\vd{q})$ is
\begin{equation}
    \vd{b}(\vd{q}_1,\vd{q}_2) = \frac{1}{2}\left[ \vd{n}(\vd{q}_1 + \vd{q}_2) - \vd{n}(\vd{q}_1) - \vd{n}(\vd{q}_2) \right] \text{.}
\end{equation}

\subsection{Affine parametric dependence} \label{sec:nonlinear:parameter}
Until now, we have suppressed the dependence of $\td{A}$ on the parameter vector $\vd{\mu} \in \mathbb{R}^{N_\mu}$. Here, we show that by precomputing certain operators, the matrices in the ROM can quickly be computed given a new parameter vector. We assume this parameter dependence is affine, i.e., $\td{A}(\vd{\mu})$ can be expressed as a linear combination of $M_\mu$ matrices as
\begin{equation}
    \td{A}(\vd{\mu}) = \sum_{j = 1}^{M_{\mu}}\zeta_j(\vd{\mu})\td{A}_j \text{,}
\end{equation}
with scalar coefficients $\zeta_j$ that depend on $\vd{\mu}$. The components of the ROM that inherit this parameter dependence are $\td{E}_k$, $\td{J}_k$, $\td{N}_k$, $\td{M}_k$, $\vd{n}_{klmn}$, $\vd{m}_{klmn}$, $\td{H}_k$, and $\tilde{\td{A}}_k$. The parameter dependence of all but the latter two is through the matrix $\td{G}^+_k$. 
\\

We define the matrices
\begin{equation}
    \td{G}_{kj} = \left( \delta_{1j} i\omega_k \td{I} - \td{A}_j \right) \tdh{Q}_k \text{,}
\end{equation}
for $j \in \{1, \dots, M_\mu \}$, where $\delta$ is the Kronocker delta. The weighted sum of these matrices is $\td{G}_k$, i.e., $\sum_{j = 1}^{M_\mu} \zeta_j (\vd{\mu}) \td{G}_{kj} = \td{G}_k$, and they are elements of $\mathbb{C}^{N_x \times N_d}$. We then define the following quantities, 
\begin{subequations}
    
\begin{equation}
    \td{G}^G_{kij} = \td{G}^*_{ki} \td{W} \td{G}_{kj} \in \mathbb{C}^{N_d \times N_d} \text{,}
\end{equation} 

\begin{equation}
    \td{B}^G_{kj} = \td{G}^*_{kj} \td{W} \td{B} \in \mathbb{C}^{N_d \times N_f}\text{,}
\end{equation}

\begin{equation}
    \td{U}^G_{kj} = \td{G}^*_{kj}\td{W} \td{U}^n \in \mathbb{C}^{N_d \times p_2}\text{,}
\end{equation}

\begin{equation}
    \vd{b}^G_{kjlmin} = \td{G}^*_{kj} \td{W} \vd{b}(\vd{\psi}_{lm},\vd{\psi}_{in}) \in \mathbb{C}^{N_d} \text{,}
\end{equation}

\begin{equation}
    \tilde{\td{A}}_{kj} = \SPOD{k}^* \td{W} \td{A}_j \SPOD{k} \text{.}
\end{equation}

\end{subequations}
Here, the superscript $`G$' and subscript $`kj'$ are reminders that the quantity comes from taking inner products with the columns of $\td{G}^*_{kj}$. With these quantities precomputed, the ROM components may be computed quickly online given a parameter vector $\vd{\mu}$ using the formulae
\begin{subequations} \label{eq:affine:ROMComps}
    
\begin{equation}
   \td{E}_k = \SPOD{k}^* \td{W} \tdh{Q}_k  \left(\td{G}_k^* \td{W} \td{G}_k\right)^{-1} \sum_{j=1}^{M_\mu} \zeta_j(\vd{\mu})\td{B}^G_{kj} \text{,}   
\end{equation}

\begin{equation}
    \td{J}_k =\td{\Phi}^*\td{W} \tdh{Q}_k  \left(\td{G}_k^* \td{W} \td{G}_k\right)^{-1} \sum_{j=1}^{M_\mu} \zeta_j(\vd{\mu})\td{B}^G_{kj} \text{,}   
\end{equation}

\begin{equation}
    \td{N}_k = \SPOD{k}^* \td{W} \tdh{Q}_k  \left(\td{G}_k^* \td{W} \td{G}_k\right)^{-1} \left( \sum_{j=1}^{M_\mu} \zeta_j(\vd{\mu})\td{U}^G_{kj} \right) \left(\td{P}^{nT} \td{U}^n \right)^{-1} \text{,}   
\end{equation}

\begin{equation}
    \td{M}_k = \td{\Phi}^* \td{W} \tdh{Q}_k  \left(\td{G}_k^* \td{W} \td{G}_k\right)^{-1} \left( \sum_{j=1}^{M_\mu} \zeta_j(\vd{\mu})\td{U}^G_{kj} \right) \left(\td{P}^{nT} \td{U}^n \right)^{-1} \text{,}   
\end{equation}

\begin{equation}
    \vd{n}_{klmn} = \SPOD{k}^* \td{W} \tdh{Q}_k  \left(\td{G}_k^* \td{W} \td{G}_k\right)^{-1} \sum_{j = 1}^{M_\mu} \zeta_j(\vd{\mu}) \vd{b}^G_{kjlmin} \text{,}
\end{equation}

\begin{equation}
    \vd{m}_{klmn} = \td{\Phi}^* \td{W} \tdh{Q}_k  \left(\td{G}_k^* \td{W} \td{G}_k\right)^{-1} \sum_{j = 1}^{M_\mu} \zeta_j(\vd{\mu}) \vd{b}^G_{kjlmin} \text{,}
\end{equation}

\begin{equation}
    \tilde{\td{A}}_k = \sum_{j = 1}^{M_\mu} \zeta_j(\vd{\mu}) \tilde{\td{A}}_{kj} \text{,}
\end{equation}

\end{subequations}
where $(\td{G}_k^* \td{W} \td{G}_k)^{-1} \in \mathbb{C}^{N_d \times N_d}$ is computed quickly using
\begin{equation}
    \left(\td{G}_k^* \td{W} \td{G}_k\right)^{-1} = \left( \sum_{j = 1}^{M_\mu} \sum_{l = 1}^{M_\mu} \zeta_j(\vd{\mu})\zeta_l(\vd{\mu}) \td{G}^G_{kjl} \right)^{-1} \text{.}
\end{equation}
Note that though $\td{H}_k$ depends on $\vd{\mu}$, this dependence is only through $\tilde{\td{A}}_k$, so after $\tilde{\td{A}}_k$ is computed with a given parameter, $\td{H}_k$ can be computed quickly using \eqref{eq:Lin:Hdef}. With the formulae given in \eqref{eq:affine:ROMComps}, all ROM components may be computed given a parameter $\vd{\mu}$ without performing operations that scale with $N_x$, as desired.

\subsection{Solving the nonlinear system} \label{sec:nonlinear:solve}
Equations \eqref{eq:nonlinear:deim:system} and \eqref{eq:nonlinear:bilinear:system}, pertaining to the DEIM-based and triadic interaction-based approximations of the nonlinearity, respectively, are both systems of $rN_\omega$ nonlinear equations in $rN_\omega$ unknowns. The online phase of the proposed reduced-order model consists of solving these equations for the coefficients $\vdt{a}_\set{K}$. For the method to be viable, these equations must be solved as quickly as possible, and to this end, we propose a fixed-point iteration technique. In the majority of our tests, the fixed-point iteration converged in a few ($\sim 10$) iterations. We also describe a slower but more stable pseudo-time-stepping method that we used in the few cases where the fixed-point iteration did not converge. We analyze the convergence of both in Appendix~\ref{app:solution_analysis}.
\\

Both systems \eqref{eq:nonlinear:deim:system} and \eqref{eq:nonlinear:bilinear:system} may be written abstractly as
\begin{equation} \label{eq:nonlinear:solve:abstract_form}
    \vdt{a}_\set{K} = \vd{c}_\set{K} + \vd{w}_\set{K}(\vdt{a}_\set{K}) \text{,}
\end{equation}
where $\vd{c}_\set{K} \in \mathbb{C}^{rN_\omega}$ is the grouping of terms in the system (either \eqref{eq:nonlinear:deim:system} or \eqref{eq:nonlinear:bilinear:system}) that do not depend on $\vdt{a}_\set{K}$, and where $\vd{w}_\set{K}: \mathbb{C}^{rN_\omega} \to \mathbb{C}^{rN_\omega}$ is the grouping of terms that do depend on $\vdt{a}_\set{K}$. The fixed-point iteration is simply 
\begin{equation} \label{eq:nonlinear:solving:iteration}
\begin{aligned}
    &\vd{a}_\set{K}^{0} = \vd{0} \text{,} \\
    &\vd{a}_\set{K}^{i+1} = \vd{c}_\set{K} + \vd{w}_\set{K}(\vd{a}_\set{K}^{i}) \text{.}
\end{aligned}
\end{equation}
With the initial guess of $\vd{0}$, the first iterate is $\vd{a}_\set{K}^{1} = \vd{c}_\set{K}$, the solution to the system without the presence of the nonlinearity. The solution $\vdt{a}_\set{K}$ to \eqref{eq:nonlinear:solve:abstract_form} is a fixed point of the iteration above.
\\

A necessary condition for convergence to this fixed point is that the eigenvalues of the Jacobian about it must be within the unit circle. In our numerical examples, we have found that this condition is usually met, and the iteration converges. The exceptions have been in cases where the nonlinearity is strong, and the solution to the system with no nonlinearity bears little resemblance to that with the nonlinearity. In Appendix~\ref{app:solution_analysis}, we present some analysis of the fixed point iteration that is consistent with this observation. We also note that Anderson acceleration \cite{Anderson65,Walker11} may be used to accelerate the convergence and possibly make it more robust.
\\

In cases where the fixed point iteration does not converge, a slower but more stable alternative is the following pseudo-time-stepping method used by Ref.~\cite{Hall02} to solve similar problems,
\begin{equation}\label{eq:nonlinear:solving:pseudotime}
\begin{aligned}
     &\vd{a}_\set{K} (0) =  \vd{c}_\set{K} \text{,} \\   
    &\frac{\text{d}}{\text{d} \tau} \vd{a}_\set{K} (\tau) = \vd{c}_\set{K} + \vd{w}_\set{K}( \vd{a}_\set{K} (\tau)) -  \vd{a}_\set{K} (\tau) \text{.}
\end{aligned}
\end{equation}
The stability of this method also depends on the eigenvalues of the same Jacobian. In this case, however, the stability condition is that the real part of the eigenvalues must be less than $1$ (with an infinitesimal time step), so the pseudo-time-stepping method is more stable than the fixed point iteration. In our numerical examples, this method never failed to converge. Note that the pseudo-time-stepping method is equivalent to the fixed point iteration if an explicit Euler integrator is used with a time step of $1$.

\subsection{Scaling analysis} \label{sec:nonlinear:scaling}
Algorithms for the offline and online phases of the method are given in Appendix~\ref{app:algorithms}. Here we give the scalings for both, depending on the handling of the nonlinearity. With either DEIM or the sparse triadic interactions, the dominant online cost comes from calculating the nonlinear term $\vd{w}_\set{K}$ at each iteration. In totaling the online scaling of the method in this subsection and timing the method in Section~\ref{sec:results}, we count all operations beginning from the initial condition and forcing and ending at the SPOD coefficients. 
\\

The online phase of the method consists of repeating the iteration \eqref{eq:nonlinear:solving:iteration} until the SPOD coefficients are converged. The constant $\vd{c}_\set{K}$ is computed once at the beginning; computing this constant is equivalent to solving for the SPOD coefficients with no nonlinearity present. The nonlinear term is computed for each iteration, and, in practice, it is the dominant cost of the method. 
\\

CPU cost due to the constant term scales as 
\begin{equation}
\mathcal{O}(N_\omega (N_f \log N_\omega + p_1N_f + rN_f + rp_1) + p_1N_x) \text{.}
\end{equation}
The first term is due to the FFT of the forcing; the second to calculating the effect of the forcing in the intermediary basis; the third to calculating the effect of the forcing directly on each SPOD coefficient; and the fourth to calculating the effect of each coefficient in the intermediary basis on each SPOD coefficient. The last term, which does not scale with the number of frequencies, comes from calculating the initial condition in the intermediary basis. 
\\

The cost of the nonlinear term at each iteration, i.e., each evaluation of $\vd{w}_\set{K}(\vdt{a}_\set{K})$, depends on the handling of the nonlinearity. If DEIM is used, each iteration scales as
\begin{equation} \label{eq:nonlinear:scaling:iter_deim}
    \mathcal{O}\left(N_\omega (r p_2 + p_2 \log N_\omega + p_1p_2 + rp_1) \right) \text{.}
\end{equation}
The first term is due to sampling the trajectory at the $p_2$ sample points at every frequency; the second to taking the IFFT of this sampling, computing the nonlinearty, then taking the FFT; the third to computing the effect of the nonlinearity at $p_2$ sample points in the intermediary basis; and the fourth to calculating the effect of the each coefficient in the intermediary basis on each SPOD coefficient. 
\\

If, instead, sparse triadic interactions are used, the scaling of each iteration is 
\begin{equation} \label{eq:nonlinear:scaling:iter_bil}
    \mathcal{O}\left( N_\omega p_1 r + \sum_{k \in \set{K}} |\set{N}_k |(r_k + p_1) \right) \text{.}
\end{equation}
The first term accounts for the effect of each coefficient of the intermediary basis on each SPOD coefficient. The sum accounts for the number of nonlinear terms that are retained at each frequency.
\\
 
The total online cost, therefore, scales as
\begin{equation}
    \mathcal{O} \left(   N_\omega \left( N_f\log N_\omega + p_1 N_f + rN_f + N_{iter} \left( rp_2 + p_2 \log N_\omega + p_1 p_2 + rp_1 \right)   \right) + p_1N_x \right)
\end{equation}
if DEIM is used and 
\begin{equation}
    \mathcal{O} \left(   N_\omega \left( N_f\log N_\omega + p_1 N_f + rN_f\right) + N_{iter}\left( rp_1N_\omega  + \sum_{k \in \set{K}} |\set{N}_k |(r_k + p_1)   \right) + p_1N_x \right)
\end{equation}
if the sparse triadic interactions are used. 
\\

In practice, in the DEIM case, $r$, $p_1$, and $p_2$ will likely be of similar magnitudes, and the dominant cost will come from the iteration to solve the system. The scaling will then be 
\begin{equation}
    \mathcal{O} \left( N_\omega N_{iter} r^2 \right) \text{.}
\end{equation}
If the sparse triadic interactions are used, the dominant cost will likely come from evaluating the nonlinear term at every iteration. A rough estimate of the scaling in this case is
\begin{equation}
    \mathcal{O}({N_{iter}N_{int}r}) \text{,}
\end{equation}
where $N_{int} = \sum_{k \in \set{K}} |\set{N}_k |$. In practice, we have observed the fixed-point iteration to converge in $\sim 10$ iterations or fewer. If the pseudo-time-stepping method is used, $N_{iter}$ is the number of time steps needed for convergence, and in the cases where we used this method, $N_{iter}$ was a few hundred using an adaptive ttime-steppingmethod (though we made no attempt to optimize this). 
\\

Assuming the governing equations are sparse, i.e., computing the action $\td{A}$ and $\vd{n}$ both scale as $\mathcal{O}(N_x)$, the offline scalings are as follows. Obtaining the SPOD modes and computing the operators necessary for the constant term scale as
\begin{equation}  
\mathcal{O}(N_\omega N_d^2 N_x) \text{.}
\end{equation}
Constructing the operators needed for the nonlinearity in the DEIM case scales as
\begin{equation}
    \mathcal{O}((r + p_1 + p_2)N_dN_xN_\omega )\text{.}
\end{equation}
In the case of sparse triadic interactions, constructing the necessary vectors scales as 
\begin{equation}
    \mathcal{O}\left(\sum_{k \in \set{K}} \sum_{l \in \set{K}} r_l r_i r_k N_x \right) \text{.}
\end{equation}
\\

So long as the offline cost is feasible, the online cost is the salient quantity. In practice, computing the effect of the nonlinearity for each iteration (\eqref{eq:nonlinear:scaling:iter_deim} or \eqref{eq:nonlinear:scaling:iter_bil}) dominates the online cost. In our numerical examples, we find that DEIM is faster (at the same accuracy) than the sparsified triadic interactions, though we do not expect this to hold in general.

\section{Results} \label{sec:results}
We use the 1-dimensional Ginzburg-Landau equations to test the method. We use the standard form of these equations for testing the DEIM-based handling of nonlinearity and a modified Ginzburg-Landau system with a quadratic nonlinearity to test the triadic-interaction-based treatment of the nonlinearity. For both systems, we compare the relative error and online CPU time of SSOP to those of POD-G. The results are encouraging: for the same number of modes, SSOP gives two orders of magnitude lower error than does POD-G. For the DEIM-based treatment of the nonlinearity, this accuracy improvement comes at slightly lower CPU time than POD-G, whereas the triadic-interaction-based ROM is slower than POD-G. We also test the accuracy of the method on out-of-sample data and find that the method trained on one GL system generalizes well to a new GL system. 
\\
\subsection{Standard Ginzburg-Landau system}
The complex Ginzburg-Landau equation is a common test case for model reduction methods \cite{Ilak10,Brunton14,Chen21,Padovan24}. It is
\begin{equation} \label{eq:results:GL_standard}
    \dot{q}(x,t) = \bigg[-\nu \frac{\partial}{\partial x} + \gamma \frac{\partial^2}{\partial x^2} + \mu_0 - c_\mu^2 + \frac{\mu_2}{2}x^2 \bigg] q(x,t) - \alpha q(x,t) |q(x,t)|^2  + f(x,t)\text{,}
\end{equation}
where $f(x,t)$ is a forcing and where we set $\nu = 2 + 0.4i$, $\gamma = 1-i$, $c_\mu = 0.2$, $\mu_2 = -0.01$, and $\alpha = 1$, which are standard values \cite{Bagheri09,Towne2018spectral}. In most of our tests, we set $\mu_0$, a bifurcation parameter in the system, to $\mu_0 = 0.229$ \cite{Bagheri09,Towne2018spectral}. The terms in the system generate advection, diffusion, and local growth or decay depending on the sign of $\mu_0 - c_\mu^2 + \frac{\mu_2}{2}x^2 $. For $\mu_0 = 0.229$, $\mu_0 - c_\mu^2 + \frac{\mu_2}{2}x^2 > 0$ for $x \in [-6.15,6.15]$, so solutions grow in this region. For this value of $\mu_0$, the system is linearly stable, i.e., the eigenvalues of the linear operator about the steady solution $q(x) = 0$ are all in the stable half-plane. There is a bifurcation at $\mu_0 = 0.397$; above this value, the system is linearly unstable, and the dynamics are characterized by a competition between the amplifying effect of the linear terms and the damping effect of the nonlinear term. While most of our tests use $\mu_0 = 0.229$, we also test the range $\mu_0 \in [0.079, 0.499]$, which envelopes the bifurcation. We refer to \eqref{eq:results:GL_standard} as the `standard' Ginzburg-Landau system to distinguish it from the modified system introduced in the next subsection.
\\

The full-order model consists of a pseudo-spectral Hermite discretization \cite{Bagheri09,Chen11} of \eqref{eq:results:GL_standard} with $N_x = 220$ collocation points \cite{Towne2018spectral} that we integrate with \texttt{ode45} in MATLAB. To generate data, we integrate for $3000$ time steps with $\Delta t = 0.8$. We use the method described in Ref. \cite{Towne2018spectral} to obtain the modes from the single long trajectory by forming $44$ trajectories of length $N_\omega = 256$, each overlapping $75 \%$ with the next. The forcing generating this data is spread over the range $x\in [-12,-8]$ and is strongest at $\overline{x}=10$. It is stochastic with a spatial correlation length of $1$, a temporal correlation length of $3.33$, and a Gaussian support. The spatiotemporal correlation is given by
\begin{equation} \label{eq:results:stoch_forcing}
    \mathbb{E}[f(x_1,t_1) f^*(x_2,t_2)] \propto \exp \left[ - \left( (x_1 - \overline{x})^2 + (x_2 - \overline{x})^2 + (x_2  - x_1)^2 + (0.3(t_2 - t_1))^2 \right)   \right] \text{,}
\end{equation}
and we take it to be zero outside $x\in [-12,-8]$.
\\

\if\showfigs1
\begin{figure}[]
    \centering
    \includegraphics{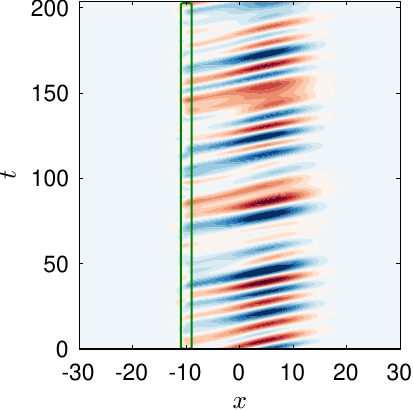}
    \caption{A trajectory of the (standard) Ginzburg-Landau system with $\mu_0 = 0.229$. The green lines demarcate the location of the forcing. The spatiotemporal structure enables space-time modes, like SPOD modes, to represent the trajectory far more efficiently than space-only modes.}
    \label{fig:state_spacetime}
\end{figure}
\fi

Figure~\ref{fig:state_spacetime} shows a space-time diagram of a trajectory of the Ginzburg-Landau system with $\mu_0 = 0.229$. The advective behavior of the system is evident in the diagonally oriented streaks, with the upward inclination indicating advection in the positive direction. There are spatial correlations here, but, crucially, there are strong spatiotemporal correlations as well. These spatiotemporal correlations --- the fact that there is structure in this space-time diagram beyond the spatial structure --- are what allow space-time modes to encode the trajectory substantially more efficiently than POD modes. 
\\

\if\showfigs1
\begin{figure}[]
    \centering
    \includegraphics{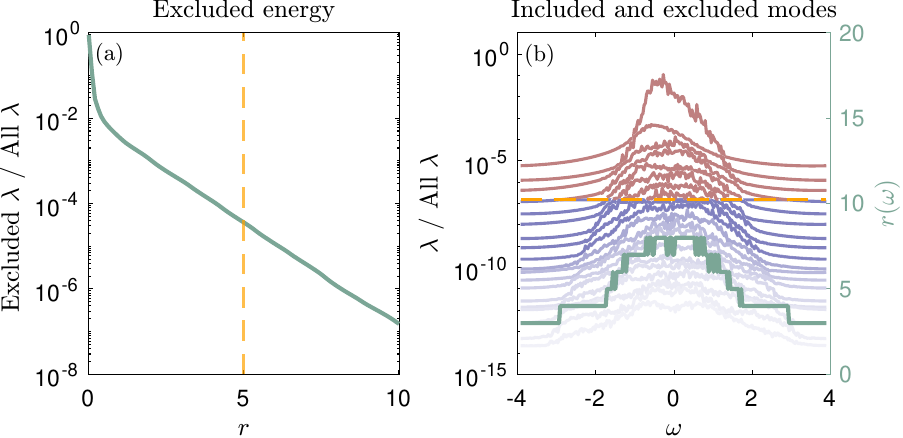}
    \caption{SPOD energies. (a) The energy not captured by the first $rN_\omega$ modes, calculated as the ratio of the sum of the energies of the excluded modes to the sum of the energies of all the modes. The dashed line corresponds to the $r$ value used to determine the cutoff in (b). (b) The energy of the included (red) and excluded (blue) modes as a function of $\omega$ for $r = 5$. Inclusion of a given mode $\vd{\psi}_{km}$ is determined by whether the associated energy is among the $rN_\omega$ largest energies over all frequencies i.e., whether $\lambda_{km} \geq \tilde{\lambda}_{(rN_\omega)}$ (see \eqref{eq:SPOD:rk}). The cutoff energy $\tilde{\lambda}_{(rN_\omega)}$ is shown in orange. The number of modes at each frequency $r(\omega)$ is shown in green on the right axis; the mean is $r = 5$.}
    \label{fig:excluded_energy}
\end{figure}
\fi

Figure~\ref{fig:excluded_energy} shows the energy of the modes. Specifically, Figure~\ref{fig:excluded_energy}~(a) shows the fraction of energy that is excluded as a function of $r$, which is computed by summing the energies of the modes that are not among the $rN_\omega$ most energetic, normalized by the sum of all of the energies. This is equivalent to the average relative SPOD reconstruction error over the training data. This gives a lower bound of the relative error of the method applied to the training data --- the error of the SSOP ROM cannot be lower than the SPOD reconstruction error. Figure~\ref{fig:excluded_energy}~(b) shows the SPOD energies of the modes as a function of frequency for $r=5$. The top curve is $\lambda_{\set{K},1}$, the first SPOD mode as a function of frequency, and the lower curves are higher mode numbers. The red and blue curves represent retained and truncated modes, respectively. The threshold is determined by $r$ as $\tilde{\lambda}_{(rN_\omega)}$ and is shown in orange. The number of modes at each frequency that meet the threshold and are thus retained is shown in green on the right axis. For the least energetic frequencies, as few as $3$ modes are retained, whereas for the most energetic frequencies, as many as $8$ are retained. The average value is $r=5$, by definition. 
\\

Throughout the remainder of this section, the (relative) errors are defined as follows. The error at time $t_j$ is the square norm of the difference between the ROM and FOM solutions normalized by the mean square norm of the FOM solution averaged over time and over test trajectories. This may be written as
\begin{equation} \label{results:err_t}
    e_j = \frac{\| \vdt{q}_j - \vd{q}_j  \|_{x}^2}{\frac{1}{N_\omega N_{test}} \sum_{i = 1}^{N_{test}} \| \vd{q}^{i}_\set{J} \|_{x,t}^2} \text{,}
\end{equation}
where $\vd{q}$ is the FOM solution, $\vdt{q}$ is the ROM approximation thereof, and $\boldsymbol{q}_\set{J}^i$ is the $i$-th trajectory in the test data. This error averaged over time is
\begin{equation}
    e = \frac{\| \vdt{q}_\set{J} - \vd{q}_\set{J}  \|_{x,t}^2}{\frac{1}{N_{test}}\sum_{i = 1}^{N_{test}} \| \vd{q}^{i}_\set{J} \|_{x,t}^2} \text{.}
\end{equation}
For the most part, we report these quantities averaged over the test trajectories below.
\\

We first test the method on $30$ new trajectories, i.e., trajectories not included in the training data. The forcing is stochastic in each case and is drawn from the same distribution as the stochastic forcing used to generate the training data (but the realizations are different). The $30$ initial conditions in the test data are taken from a single long run of the system and are also not in the data used to compute the SPOD modes, but come from the same statistical distribution as the states in the training data. These conditions together mean that the SPOD modes from the training data are good at representing the trajectories in the test data.
\\

\begin{figure}[!h]
    \centering
    \includegraphics{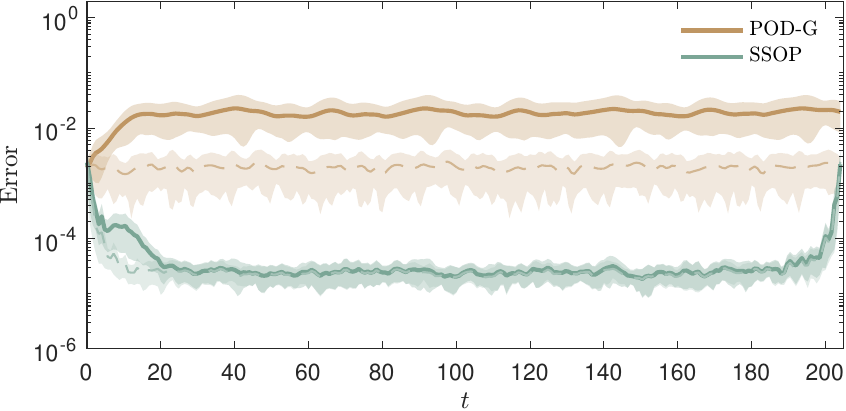}
    \caption{Error as a function of time for the proposed method and for POD-G, both with $r=5$. The curves represent the mean error over the $30$ test trajectories, and the shaded regions indicate the standard deviation of these errors. The dashed curves show the respective projection errors, which serve as lower bounds for the two methods.}
    \label{fig:cubic_errtime}
\end{figure}

In Figure~\ref{fig:cubic_errtime}, we show the relative error of the proposed method with $r=5$ modes along with that of POD-G, also with $5$ modes. The dashed curves are the projection errors of the respective bases. In other words, they represent the discrepancy between the full-order solution and the projection thereof onto the POD and SPOD bases. The projection error is a lower bound for the model error; if the model recovers the exact coefficients, then it achieves the projection error. All errors are calculated as the mean over the $30$ realizations on which we test the method, and the shaded regions represent data within one standard deviation of the mean. Most notably, SSOP outperforms POD-G by nearly three orders of magnitude for most of the interval. It also significantly outperforms the projection of the FOM solution onto the POD modes (dashed brown). Since POD modes are the optimal linear encoder, the POD projection error is a lower bound for any method based on a space-only linear encoding, such as a space-only Petrov-Galerkin method or operator inference. That SSOP and the SPOD projection (solid and dashed green, respectively) nearly overlap indicates that SSOP solves nearly exactly for the SPOD coefficients. A notable feature of Figure~\ref{fig:cubic_errtime} is that the error is largest at the ends of the interval, and it bears mentioning that this feature is not due to Gibbs phenomena --- indeed, it is too small. Instead, it is due to a cyclic property of the SPOD encoder/decoder pair, and as more modes are included, the effect diminishes. 
\\

\begin{figure}[!h]
    \centering
    \includegraphics{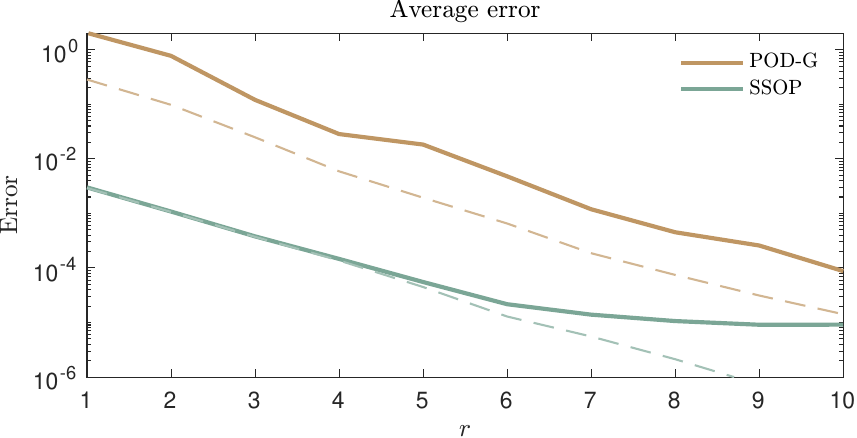}
    \caption{Error, averaged over time, for SSOP and for POD-G as a function of $r$, the number of modes used in the ROMs. The dashed curves again represent the projection error onto the POD and SPOD modes.}
    \label{fig:cubic_err}
\end{figure}

Figure~\ref{fig:cubic_err} shows the relative error, averaged over time, as a function of the number of modes retained. For all the values shown, SSOP again gives multiple orders of magnitude lower error than POD-G. The dashed curves are again the projection error. Another notable feature of the method is that it nearly achieves the SPOD mode projection error for relatively few modes. The error asymptote is likely due to the fact that the forcing and nonlinearity are not exactly represented by a finite Fourier series, as was assumed in deriving the method. 
\\

\begin{figure}[!h]
    \centering
    \includegraphics{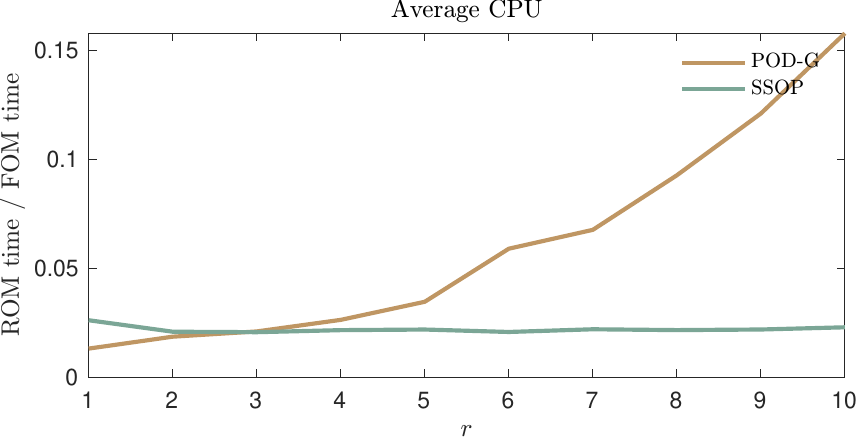}
    \caption{CPU time for SSOP in comparison to that of POD-G as a function of the number of modes used in the ROMs.}
    \label{fig:cubic_scaling}
\end{figure}

Figure~\ref{fig:cubic_scaling} shows the CPU time of the proposed method in comparison to that of POD-G. The substantial error reduction shown does not come at an increase in computational cost; indeed, the method is faster for most of the range considered. In what follows, we use $r=5$ for all tests. 
\\

Next, we investigate how the method performs for non-stationary forcings.  We again use the SPOD modes from data generated by a stochastic, statistically stationary forcing, but test the resulting model on a variety of non-stationary forcings. The forcing is active in the same spatial region as before, and is given in this region by
\begin{equation}
    f(x,t) = f(t) \sqrt{\exp[-(x - \overline{x})^2]} \text{,}
\end{equation}
The forcing is again strongest at $\overline{x} = 10$. We test four choices of the function $f(t)$: a periodic function, a pulse function, a quasiperiodic function, and a series of pulses, steps, and quasiperiodic functions. The initial condition is zero in each case. 
\\

\begin{figure}[]
    \centering
    \includegraphics{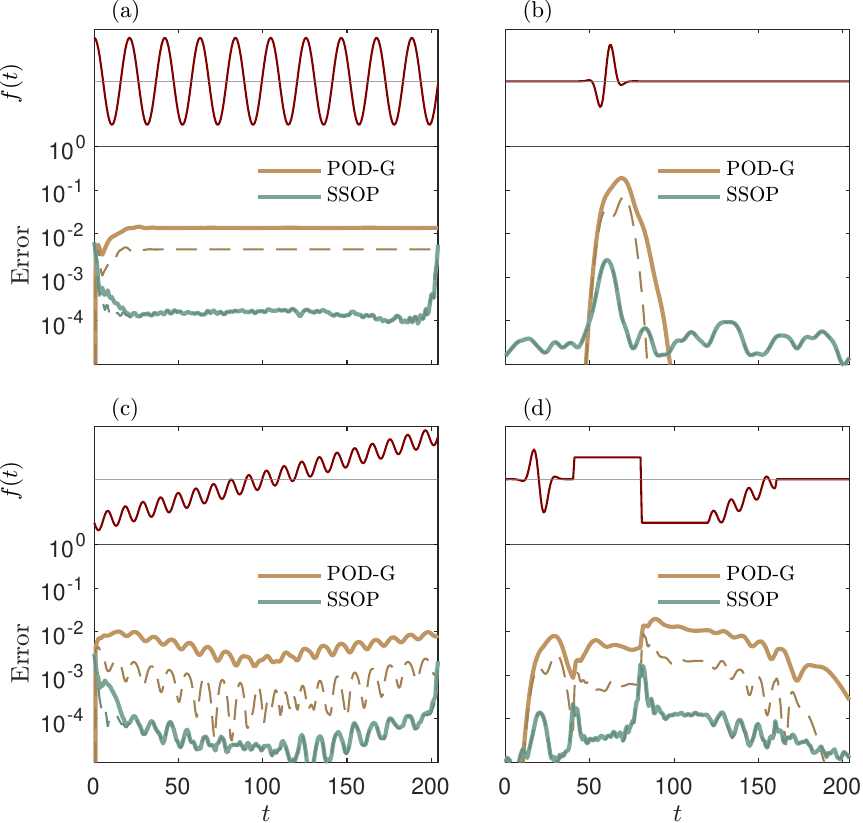}
    \caption{The error as a function of time for a variety of non-stochastic forcings for both SSOP and POD-G with $r=5$ modes, along with the respective projection errors (dashed). The SSOP performance is robust for these out-of-sample forcings.}
    \label{fig:diff_forcings}
\end{figure}

These functions, the relative error of SSOP and POG-G, and the respective projection errors are shown in Figure~\ref{fig:diff_forcings}. The salient takeaways are the same as in the stochastic forcing case: SSOP achieves substantially lower error than POD-G and than the POD projection error, and nearly achieves the SPOD projection error. With the periodic (a) and quasiperiodic (c) forcings, the solution at the beginning of the interval is different than the solution at the end of the interval, so the SPOD projection error, and thus SSOP, is high at the beginning and end of the interval. Conversely, because the pulse (b) and series (d) forcing are both zero for long enough at the end of the interval for the solution to decay, the error of the SPOD projection and SSOP remains low at the beginning and end of the interval. Note that because the error is normalized by the average square norm of the solution over the interval, the pulse error is substantially higher than the others. We also note that none of the forcings here (or in the stochastic forcing case) meet assumption \textit{(ii)}, but this has little effect on the accuracy of the method. 
\\

Finally, we test the method on a range of Ginzburg-Landau systems by varying $\mu_0$ over the range $\mu_0 \in [0.079,0.499]$ in increments of $0.03$. The behavior of the Ginzburg-Landau system changes significantly over this range. At the lower end, the system is linearly stable, i.e., all of the eigenvalues of $\td{A}$ are in the stable half plane, and the behavior is entirely modal, i.e., the eigenvectors are nearly orthogonal. At $\mu_0 = 0.229$, the value we have used in the tests so far, the system is slightly non-modal --- one measure of this is the optimal transient growth \cite{Trefethen93,Schmid07}, which is nearly $5$. At $\mu_0 = 0.379$, the system is linearly stable but is strongly non-modal, with an optimal transient growth of nearly $200$. At $\mu_0 > 0.397$, the system is linearly unstable, and the dynamics feature a competition between the growth caused by the linear instability and the damping of the nonlinear term. 
\\

\begin{figure}[]
    \centering
    \includegraphics{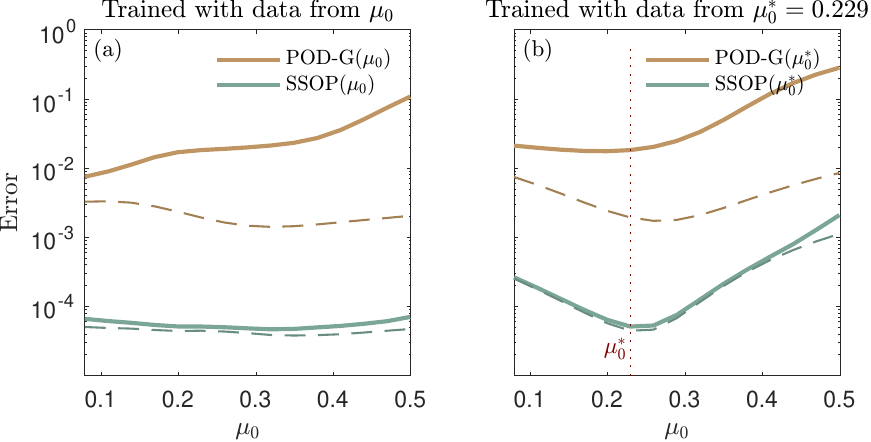}
    \caption{The methods tested on a range of Ginzburg-Landau systems parameterized by $\mu_0$ with $r=5$. In (a), the model is trained using data from the same GL system. In (b), the models are trained using data from $\mu_0^* = 0.229$. The dashed curves correspond to the projection error for the modes from the training system.}
    \label{fig:vary_mu}
\end{figure}

We perform two tests over this range of $\mu_0$. First, we build an SSOP model with SPOD modes from training data at each $\mu_0$ with the stochastic forcing \eqref{eq:results:stoch_forcing}, then test the model at the same $\mu_0$ with different realizations from the same distribution of forcings. This test verifies that the method can perform well across a range of Ginzburg-Landau systems. Second, we take SPOD modes from $\mu_0^* = 0.229$, then use them to build an SSOP model for each $\mu_0$. This second test is designed to see if the method can be used to take data from one system and use it to predict behavior of a different one.
\\

Figure~\ref{fig:vary_mu}(a) shows the results of the first test using $r=5$ modes --- the average relative error over both time and the $30$ trajectories is shown as a function of $\mu_0$. Again, SSOP substantially outperforms POD-G and the projection error of the POD basis, and nearly achieves its own projection error. The SSOP error is remarkably uniform over the range of parameters, in contrast to the POD-G error. The SSOP performance past the bifurcation at $\mu_0 = 0.397$ shows the efficacy of the strategy described in Section~\ref{sec:nonlinear:implicit} of diverting some of the linear term to the nonlinear term. 
\\

Figure~\ref{fig:vary_mu}(b) shows the second test with $r= 5$ modes. Again, the relative error averaged over time and trajectories is shown, but, in contrast to (a), the model is built using modes educed using data from a system with $\mu_0^* = 0.229$. The projection errors (dashed) are also computed using the modes from $\mu_0^*$. That the projection error is substantially lower near $\mu_0^*$ than at other values indicates that the SPOD modes change significantly as $\mu_0$ changes. Despite this, the model error results (solid) show SSOP can generate accurate predictions for new systems. In fact, the SSOP model built using data from a different Ginzburg-Landau system ((b), solid green) produces lower error in this example than the projection error onto the POD modes of the test system ((a), dashed brown). We view this as a strong result showcasing the robustness of the proposed method.

\subsection{Quadratic Ginzburg-Landau system} \label{sec:results:quadratic}
Here, we test the method described in Subsection~\ref{sec:nonlinear:bil} for handling quadratic nonlinearities. We do this by replacing the standard nonlinear term with a Burgers’-type (quadratic) nonlinearity and find that most of the triadic interactions have a negligible impact. We exclude roughly $99\%$ of the interactions, which leads to substantial speedup with no meaningful increase in the error relative to retaining all of them. Even with excluding these interactions, however, the method is slower than using DEIM on the same problem. Our outlook is that using the triadic-interaction-sum method of handling nonlinearities will be more effective than DEIM for more complex problems based on the findings in Ref. \cite{Sipp20}.
\\

The quadratically nonlinear Ginzburg-Landau equation is
\begin{equation}
\frac{\partial}{\partial t} q(x,t) = \bigg[-\nu \frac{\partial}{\partial x} + \gamma \frac{\partial^2}{\partial x^2} + \mu_0 - c_\mu^2 + \frac{\mu_2}{2}x^2 \bigg] q(x,t) + \kappa q(x,t) \frac{\partial}{\partial x} q(x,t) \text{.}
\end{equation}
We set $\kappa = 5$ and otherwise use the same parameters from before (including $\mu_0 = 0.229$).
\\

We generate data following the same procedure as before --- using the stochastic forcing \eqref{eq:results:stoch_forcing} to generate $3000$ time steps of training data with $\Delta t  = 0.8$, forming this data into $44$ overlapping blocks of length $N_\omega = 256$, and computing the SPOD modes and operators. We then again use different realizations from the same forcing distribution to generate the data used for testing. 
\\

\begin{figure}[]
    \centering
    \includegraphics{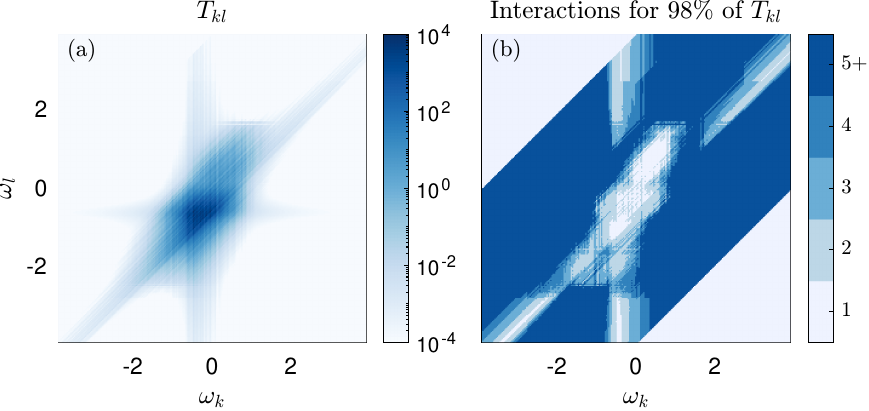}
    \caption{The impact of $\omega_l$ on $\omega_k$ as defined by \eqref{eq:results:impact} via triadic interaction with $\omega_i = \omega_k - \omega_l$. (a) shows $T_{kl}$, a proxy for the strength of all interactions between the SPOD modes at $\omega_l$ with those at $\omega_i$. (b) shows how many of these interactions are needed to account for $98 \%$ of $T_{kl}$.}
    \label{fig:int_map}
\end{figure}

The criterion (described in Section~\ref{sec:nonlinear:bil}) for including the triadic interaction between mode $m$ at $\omega_l$ and mode $n$ at frequency $\omega_i = \omega_k - \omega_l$ is $\|\vd{n}_{klmn}\|_{x}^2 \lambda_{lm}\lambda_{in} \geq \epsilon$ for a user-defined value of $\epsilon$. We define the matrix $\td{T} \in \mathbb{R}^{N_\omega \times N_\omega}$ as
\begin{equation} \label{eq:results:impact}
    T_{kl} = \sum_{m = 1}^{r_l} \sum_{n = 1}^{r_i} \|\vd{n}_{klmn}\|_{x}^2 \lambda_{lm}\lambda_{in} \text{.}
\end{equation}
$T_{kl}$ is a proxy for the total impact of $\omega_l$ on $\omega_k$ via triadic interactions. We plot $ T_{kl}$ for all pairs of frequencies in Figure~\ref{fig:int_map}(a). In Figure~\ref{fig:int_map}(b), we show the number of interactions at a given frequency pair $\omega_k$, $\omega_l$ that is required to account for $98 \%$ of $T_{kl}$. We see that for many of the frequency pairs, the total impact at a frequency pair is dominated by just a few mode pairs within the two frequencies. Noting the overlap in the high-value region of (a) and the low-value region of (b), we see that this is particularly true at the more energetic mode pairs. The implication of (a) and (b) together is that the great majority of frequencies may be discarded without having a meaningful effect on the approximation of the nonlinearity.
\\
\\

\begin{figure}[]
    \centering
    \includegraphics{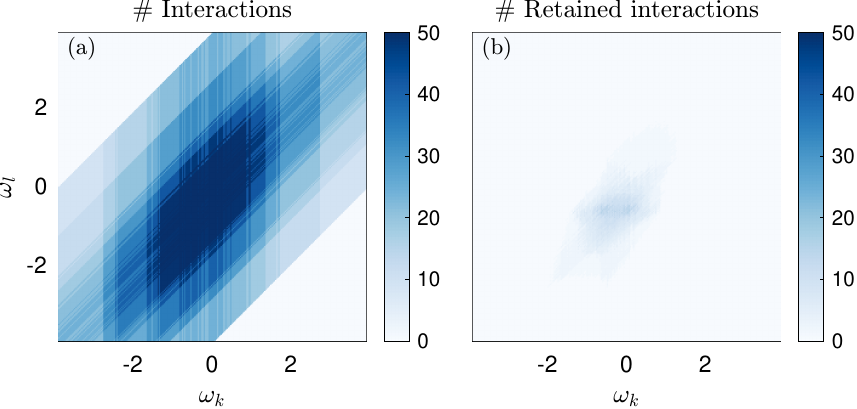}
    \caption{(a) The total number of interactions between SPOD modes at $\omega_l$ and $\omega_k - \omega_l$ with $r = 5$. (b) The number of these interactions that are retained when $\epsilon = 10^{-1.8}$. Only $1.7 \%$ of the interactions in (a) are retained with this threshold.}
    \label{fig:int_map_number}
\end{figure}

Figure~\ref{fig:int_map_number}(a) shows the number of interactions at each frequency pair with $r = 5$. These values are only a function of $r_\set{K}$ --- the number at $(\omega_k,\omega_l)$ is $r_l r_i$, where $\omega_i = \omega_k - \omega_l$. Figure~\ref{fig:int_map_number}(b) shows the number at each pair that are retained when the threshold in \eqref{eq:nonlinear:bilinear:inclusion_criterion} is set to $\epsilon = 10^{-1.8}$ (we show this one because it leads to a good balance between accuracy and speed in the results). With this choice of $\epsilon$, $1.7 \%$ of the interactions are retained. 
\\

\begin{figure}[]
    \centering
    \includegraphics{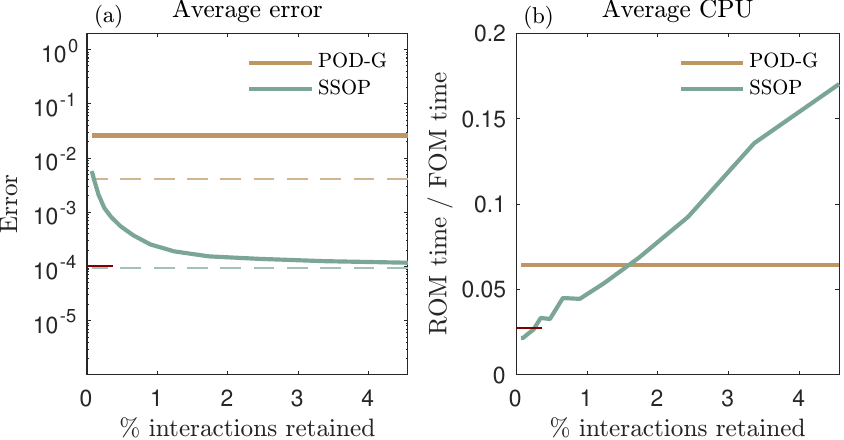}
    \caption{The relative error (a) and CPU time (b) for the method with $r = 5$ using the sparse-triadic-interaction approximation of the nonlinearity as a function of the number of triadic interactions retained. The $x$-axis shows the percentage of the total triadic interactions used. The red ticks indicate results using the DEIM approximation of the nonlinearity.}
    \label{fig:err_cpu_bil}
\end{figure}

We test the accuracy and speed of the method on the $30$ trajectories as a function of the number of interactions included. In Figure~\ref{fig:err_cpu_bil}(a), we plot the relative error averaged over time and trajectories against the percentage of the total ($1384103$) interactions. As one would expect, increasing the number of retained interactions causes the error to decrease and the CPU time to increase. Figure~\ref{fig:err_cpu_bil}(b) shows that the CPU increase is linear in the number of interactions retained, as is predicted by the scaling analysis in Section~\ref{sec:nonlinear:scaling}. The error quickly reaches a plateau as the newly retained interactions become less energetic (the plateau is above the SPOD projection error, though it is quite close). 
\\

The red ticks in Figure~\ref{fig:err_cpu_bil} indicate the error and timing results using the DEIM-based approximation of the nonlinearity. For this problem, DEIM is able to achieve lower error at lower CPU time than the triadic-interaction-based approximation, but we do not expect that this is true in general.

\section{Conclusions} \label{sec:conclusions}
We have generalized spectral solution operator projection \cite{Frame25} to nonlinear systems. The approach represents the unknown trajectory by representing each temporal frequency with the basis of SPOD modes at that frequency. By projecting an implicit solution to the nonlinear system onto a set of retained SPOD modes, we arrive at a set of nonlinear algebraic equations for the SPOD coefficients. The online phase of the method comprises solving this set of equations given the initial condition, forcing, and system parameters, which we do using a fixed-point iteration. 
\\

The method performed well in a series of tests using a 1-dimensional PDE. The key result is that, for the same number of modes, the error is two orders of magnitude lower than that of POD-G and is substantially lower than the POD projection error, which is a lower bound for time-domain Petrov-Galerkin methods. Returning to the two questions outlined in the introduction, our method achieves an affirmative answer to both. Specifically, (1) it does recover the encoding of the trajectory accurately, i.e., the SSOP error is near the SPOD projection error; (2) it does so in a similar (and often shorter) time than POD-G. We found these results to be robust -- we tested SSOP against POD-G on out-of-sample forcings, as well as parameter variations, and its error was consistently much lower than that of POD-G. The accuracy improvement does not come at an increased computational cost: the method is faster than POD-G for most numbers of modes. 
\\

Three downsides of the method bear mentioning. First, the initial step of the method is to take a Fourier transform of the forcing. This means that if the forcing is not known on the entire interval before starting the method, the method cannot be applied. Second, obtaining the SPOD modes requires more training data than obtaining, e.g., the POD modes that form the basis for many space-only methods. Finally, the method applies to a particular temporal window $[0,T]$, which is set by the block length used in calculating the SPOD modes. If predictions on a longer window are desired, the method may be repeated using the final (or near the final, as the error increases at the end of the interval) state as the new initial condition. This is somewhat cumbersome in comparison to space-only methods, which preserve the Markovian property of the governing equations. On balance, however, we are quite encouraged by the results and believe they give reason for further interest in this method and space-time ROMs more broadly.
\\

Two important questions remain: (1) For which systems do the solution methods converge? In the appendix, we show that when the effect of nonlinearity is small compared to the linear solution, the fixed-point iteration and pseudo-time-stepping method converge, and that in the opposite limit, they do not for a broad class of nonlinearities. We are optimistic about the convergence for systems in which linear methods can be used to predict nonlinear behavior, as is the case in many fluid dynamic problems. (2) How accurate will the operator approximations be for larger systems, and to what extent will a lack of accuracy here impact the error in the solution? In the linear case, we found that we were able to achieve good accuracy in a $2$-dimensional advection-diffusion problem. This problem does not have strong modal behavior, which likely makes these approximations less accurate, so we expect they will hold sufficiently in larger problems, especially ones with modal behavior. 
\\

\section{Acknowledgments} We gratefully acknowledge funding from the National Science Foundation grant No. 2237537.
\begin{appendices} 
\section{Derivation of full-order equations in the frequency domain} \label{app:full_order_fd}
Here, we derive \eqref{eq:NL_method:qhat(f,q0)} from
\begin{equation} \label{eq:app1:DFT_def}
    \vdh{q}_k = \sum_{j = 0}^{N_\omega - 1} \left( e^{\td{A}j\Delta t}\vd{q}_0 + \int_0^{j \Delta t} e^{\td{A} (j \Delta t - t')}  \td{B} \vd{f}(t') + \vd{n}(\vd{q}(t')) \ dt' \right)  e^{-i\omega_k j \Delta t} \text{.}
\end{equation}
We assume \textit{(i)} the forcing is of the form $\vd{f}(t) = \sum_{l = 0}^{N_\omega-1} \vdh{f}_l e^{i\omega_lt}$, \textit{(ii)}  that $i\omega_l \td{I} - \td{A}$ is invertible for all included $l$, and  \textit{(iii)} that the nonlinearity is of the form $\vd{n}(\vd{q} (t)) = \sum_{l = 0}^{N_\omega} \vdh{n}_l e^{i\omega_lt}$. We refer to the term involving the initial condition as $\vdh{q}_{k,ic}$ and the term involving the integral as $\vdh{q}_{k,int}$. We start by evaluating
\begin{equation}
    \vdh{q}_{k,ic} = \sum_{j = 0}^{N_\omega - 1} e^{(\td{A} - i\omega_k \td{I})j \Delta t} \vd{q}_0 \text{.}
\end{equation}
This is a matrix geometric sum, i.e., each matrix in the sum is the previous one multiplied by $e^{(\td{A} - i\omega_k \td{I}) \Delta t}$. The solution to the geometric sum is $(\td{I} - e^{(\td{A} - i\omega_k \td{I}) \Delta t})^{-1}(\td{I} - e^{(\td{A} - i\omega_k \td{I}) N_\omega \Delta t}) \vd{q}_0$. Since $e^{i N_\omega \Delta t} = 1$, this is
\begin{equation} \label{eq:FOF:qic}
   \vdh{q}_{k,ic} = \left(\td{I} - e^{(\td{A} - i\omega_k \td{I}) \Delta t}\right )^{-1} \left(\td{I} - e^{\td{A} T}\right) \vd{q}_0 \text{.}
\end{equation}
Next, we evaluate $\vdh{q}_{k,int}$. Using assumptions \textit{(i)} and \textit{(iii)}, this may be rewritten as
\begin{equation} \label{eq:app1:qkint_intermediate1}
    \vdh{q}_{k,int} = \frac{1}{N_\omega} \sum_{j = 0}^{N_\omega - 1} e^{-i\omega_k j \Delta t} e^{\td{A}j \Delta t} \int_0^{j \Delta t} \sum_{l = 0}^{N_\omega - 1} e^{(i\omega_l \td{I} - \td{A} )t'} \left( \td{B}\vdh{f}_l + \vdh{n}_l \right) \ dt' \text{.}
\end{equation}
Integrating brings out a resolvent operator and \eqref{eq:app1:qkint_intermediate1} becomes
\begin{equation}
   \vdh{q}_{k,int} = \frac{1}{N_\omega} \sum_{j = 0}^{N_\omega - 1} e^{-i\omega_k j \Delta t} \sum_{l = 0}^{N_\omega - 1} \td{R}_l \left( e^{i\omega_l j \Delta t} - e^{\td{A} j \Delta t}\right) \left( \td{B} \vdh{f}_l + \vdh{n}_l \right) \text{.}
\end{equation}
Simplifying, we have 
\begin{equation}
    \vdh{q}_{k,int} = \frac{1}{N_\omega} \sum_{j = 0}^{N_\omega - 1} \sum_{l = 0}^{N_\omega -1} \td{R}_l \left(e^{i (\omega_l - \omega_k)j \Delta t} - e^{(\td{A} - i\omega_k \td{I} )j\Delta t} \right) \left( \td{B} \vdh{f}_l + \vdh{n}_l \right) \text{.}
\end{equation}
The frequency difference term evaluates to zero for $\omega_l \neq \omega_k$, and the other term may be evaluated by noting that it is a geometric sum. With this, we have
\begin{equation}
    \vdh{q}_{k,int} = \td{R}_k \left( \td{B} \vdh{f}_k + \vdh{n}_k \right) + \frac{1}{N_\omega} \left(\td{I} - e^{(\td{A} - i\omega_k \td{I}) \Delta t}\right )^{-1} \left(\td{I} - e^{\td{A} T}\right) \sum_{l = 0}^{N_\omega -1} \td{R}_l \left( \td{B}\vdh{f}_l + \vdh{n}_l \right) \text{.} 
\end{equation}
Adding $\vdh{q}_{k,ic}$ and $\vdh{q}_{k,int}$, we recover \eqref{eq:NL_method:qhat(f,q0)}.

\section{Analysis of the solution procedure} \label{app:solution_analysis}
Here, we analyze two limiting cases of the fixed-point iteration used to find the solution of the nonlinear system of equations. While we do not have practical necessary and sufficient conditions for the convergence of the iteration, these limiting cases are consistent with our observation that the cases where the iteration did not converge were ones in which the nonlinearity was strong compared to the other terms. 

The nonlinear system of equations to be solved is given in \eqref{eq:nonlinear:solve:abstract_form}, and, in order to avoid the notational burden in this section, we rewrite it as
\begin{equation} \label{eq:app:sol:system_general}
    \vd{x} = \vd{c} + \vd{f}(\vd{x})\text{,}
\end{equation}
where $\vd{c} \in \mathbb{C}^N$ is a known constant vector, $\vd{f}: \mathbb{C}^N \to \mathbb{C}^N$ is a known nonlinear function, and $\vd{x} \in \mathbb{C}^N$ is the unknown vector to be obtained.

\subsection{Small nonlinearity}
Here, we analyze the case where the nonlinear term is small, and show a correspondence between the $i$-th fixed-point iterate and the $i$-th term of a perturbation series about $\vd{c}$. We require that the nonlinear term be analytic, so the equations to solve are written 
\begin{equation} \label{eq:app:sol:system_analytic}
    \vd{x} = \vd{c} + \epsilon \sum_{k=2}^\infty \vd{f}_k \overbrace{(\vd{x},\dots, \vd{x})}^{ k \text{ arguments}} \text{,}
\end{equation}
where $\epsilon \in \mathbb{R}$, and $\vd{f}_k$ is linear in each of its $k$ arguments. We have excluded the $k=0$ and $k=1$ terms of the Taylor series for the nonlinearity since we have assumed that constant and linear terms of the right-hand side in the governing equations are accounted for by $\vd{c}$, though the following holds if these terms are included in the sum as well.

\subsubsection{Perturbation series}
If $\epsilon = 0$, then $\vd{x} = \vd{c}$. We use the regular perturbation ansatz $\vd{x} = \vd{x}_0 + \epsilon \vd{x}_1 + \epsilon^2 \vd{x}_2 + \dots$ to search for solutions near $\vd{c}$ when $\epsilon$ is small. To find a recurrence relation for $\vd{x}_i$, we insert this perturbation ansatz into \eqref{eq:app:sol:system_analytic}, giving
\begin{equation}
    \sum_{i=0}^\infty \epsilon^i \vd{x}_i = \vd{c} + \epsilon \sum_{k=2} \vd{f}_k \left(\sum_{i_1 = 0}^\infty \epsilon^{i_1} \vd{x}_{i_1} , \dots, \sum_{i_k = 0}^\infty \epsilon^{i_k} \vd{x}_{i_k} \right). 
\end{equation}
Using the multilinearity of $\vd{f}_k$ and grouping powers of $\epsilon$, this may be rewritten as 
\begin{equation}
    \sum_{i=0}^\infty \epsilon^i \vd{x}_i = \vd{c} + \sum_{i=1}^\infty \epsilon^i \sum_{k=2}^\infty \sum_{i_1 + \dots + i_k = i-1 } \vd{f}(\vd{x}_{i_1}, \dots , \vd{x}_{i_k}). 
\end{equation}
Here, the indices in the inner sum are non-negative integers. Finally, equating like powers of $\epsilon$ on the right and left, we have the following recurrence relation
\begin{subequations}
\begin{equation}\label{eq:app:sol:analytic_recurrence_base}
    \vd{x}_0 = \vd{c},
\end{equation}
\begin{equation} \label{eq:app:sol:analytic_recurrence}
    \vd{x}_i = \sum_{k=2}^\infty \sum_{i_1 + \dots + i_k = i-1 } \vd{f}_k(\vd{x}_{i_1}, \dots , \vd{x}_{i_k}). 
\end{equation}
\end{subequations}
\eqref{eq:app:sol:analytic_recurrence} is a recurrence relation since the highest order term on the right is $\vd{x}_{i-1}$. 

\subsubsection{Fixed point iteration}
The fixed-point iteration with the analytic right-hand side is
\begin{subequations}
\begin{equation} \label{eq:app:sol:iteration_base}
    \vd{x}^0 = \vd{c}
\end{equation}
\begin{equation}
    \vd{x}^{i+1} = \vd{c} + \epsilon  \sum_{k=2}^{\infty}\vd{f}_k(\vd{x}^i,\dots, \vd{x}^i) \text{.}
\end{equation}
\end{subequations}

\subsubsection{Relation between the two}
\subsubsection*{Statement}
The $i$-th fixed point iterate is equal to the perturbation series truncated at the $i$-th order \textit{to $\mathcal{O}(\epsilon^i)$}. That is, 
\begin{equation} \label{eq:app:sol:theorem}
    \vd{x}^i = \sum_{j = 0}^i \epsilon^j \vd{x}_j + \epsilon^{i+1}\vd{r}_i \quad \text{for} \quad i \in \{0, 1, \dots \}\text{,}
\end{equation}
where $\vd{r}_i$ is an error term, and the prefactor of $\epsilon^{i+1}$ indicates that this error term is of order $\epsilon^{i+1}$ or higher. 

\subsubsection*{Proof}
We prove \eqref{eq:app:sol:theorem} inductively. The base case is trivial: the zeroth iterate is exactly equal to the zeroth term in the perturbation series, since both are equal to $\vd{c}$ (see \eqref{eq:app:sol:analytic_recurrence_base} and \eqref{eq:app:sol:iteration_base}).  We assume the statement is true for $i = n$, and show that this implies it for $i=n+1$. Namely, we assume
\begin{equation}
    \vd{x}^n = \sum_{j= 0}^n \epsilon^{j}\vd{x}_{j} + \epsilon^{n+1}\vd{r}_n
\end{equation}
With this assumption, the $n+1$-st fixed point iterate is
\begin{equation} \label{eq:analytic_long1}
    \vd{x}^{n+1} = \vd{c} + \epsilon \sum_{k=2}^{\infty}\vd{f}_k \left(\sum_{i_1 = 0}^n \epsilon^{i_1}\vd{x}_{i_1} + \epsilon^{n+1}\vd{r}_n,\dots, \sum_{i_k = 0}^n \epsilon^{i_k}\vd{x}_{i_k} + \epsilon^{n+1}\vd{r}_n \right).
\end{equation}
All terms involving $\vd{r}_n$ are of order $\epsilon^{n+2}$ or higher, so they may be grouped into the $n+1$-th error term. By grouping terms at the same order in $\epsilon$, and combining all terms of order higher than $\epsilon^{n+1}$ into the error term, \eqref{eq:analytic_long1} may be rewritten as
\begin{equation}
    \vd{x}^{n+1} = \vd{c} + \sum_{j=1}^{n+1} \epsilon^{j} \sum_{k=2}^\infty \sum_{i_1 + \dots + i_k = j-1} \vd{f}_k \left (\vd{x}_{i_1}, \dots,\vd{x}_{i_k} \right) + \epsilon^{n+2}\vd{r}_{n+1}. 
\end{equation}
Since $\vd{x}_0 = \vd{c}$ and since the inner two sums are $\vd{x}_j$ by \eqref{eq:app:sol:analytic_recurrence}, the above equation may be rewritten as
\begin{equation}
    \vd{x}^{n+1} = \sum_{j=0}^{n+1} \epsilon^j \vd{x}_j + \epsilon^{n+2}\vd{r}_{n+1}. 
\end{equation}
This completes the proof.

\subsection{Bilinear right-hand side}
Here we analyze the iteration in the case where $\vd{c} = \vd{0}$ and where $\vd{f}(\vd{x}) = \vd{b}(\vd{x}, \vd{x})$ where $\vd{b}: \mathbb{C}^N \times \mathbb{C}^N \to \mathbb{C}^N$ is bilinear. This case is relevant where the nonlinear term is quadratic and dominant, and this case is sometimes discussed in the context of resolvent analysis (see Fig. 1 of Ref. \cite{McKeon17}, a review of the subject). We show that neither the fixed-point iteration nor the pseudo-time-stepping method can converge to nonzero fixed points due to a linear instability about these points. We first analyze the fixed-point iteration 
\begin{equation}
    \vd{x}^{i+1} = \vd{b}(\vd{x}^i, \vd{x}^i) \text{}
\end{equation}
near a solution $\overline{\vd{x}}$ satisfying $\overline{\vd{x}} = \vd{b}(\overline{\vd{x}},\overline{\vd{x}})$. Writing $\vd{x}^i = \overline{\vd{x}} + \epsilon \vd{y}^i$ and inserting this into the iteration, we have
\begin{equation}
    \vd{y}^{i+1} = \epsilon \left( \vd{b}(\overline{\vd{x}},\vd{y}^i) + \vd{b}(\vd{y}^i,\overline{\vd{x}})    \right) + \epsilon^2 \vd{b}(\vd{y}^i, \vd{y}^i) \text{.}
\end{equation}
Sufficiently close to the fixed point (i.e., with $\epsilon$ sufficiently small), the dynamics are then given by 
\begin{equation}
    \vd{y}^{i+1} = \td{L}\vd{y^i} \text{,}
\end{equation}
where $\td{L}\vd{y^i} = \vd{b}(\overline{\vd{x}},\vd{y}^i) + \vd{b}(\vd{y}^i,\overline{\vd{x}}) $ is the Jabobian about the fixed point $\overline{\vd{x}}$. If all the eigenvalues of $\td{L}$ are within the unit circle, $\|\vd{y}\|$ decays as the iteration proceeds, and the fixed point $\overline{\vd{x}}$ is stable; otherwise, it is unstable. 
\\

If $\overline{\vd{x}} = \vd{0}$, then the fixed point is stable since $\td{L}$ is the zero matrix. However, if $\overline{\vd{x}} \neq \vd{0}$, the fixed point is unstable because there is an eigenvalue of $2$ associated with the eigenvector $\vd{y} = \overline{\vd{x}}$. This is true because $\td{L} \overline{\vd{x}} = \vd{b}(\overline{\vd{x}},\overline{\vd{x}}) + \vd{b}(\overline{\vd{x}},\overline{\vd{x}}) = 2\overline{\vd{x}}$. This argument can be generalized to the case when the constant is zero and the nonlinearity is any $k$-linear function. In this case, the origin is again always a stable fixed point (for $k>1$), and the linear operator about any other fixed point admits an eigenvalue of $k$ associated with the eigenvector $\overline{\vd{x}}$. 
\\

The above analysis implies that the pseudo-time-stepping method also cannot converge to non-zero fixed points in the case of $\vd{c}= \vd{0}$ and a $k$-linear right-hand side since the stability in that case requires that the eigenvalues of the same Jacobian have a real part less than $1$.

\section{Simulation-based approximation of operators}
\label{app:approxH}
In the main text, we approximated the operator $\SPOD{k}^* \td{W} \left(\td{I} -e^{(\td{A} - i \omega_k \td{I})\Delta t} \right)^{-1} \left(\td{I} - e^{\td{A}T} \right) \td{\Phi} \in \mathbb{C}^{r_k \times p_1}$ by using the SPOD modes to form a Galerkin-type approximation of $\td{A}$. Here, we detail a means of forming this operator exactly (up to FOM errors) by performing $p_1$ runs of length $T$ of the linearized FOM. The benefit of the approach is that this accounts more accurately for the transient term. The drawbacks are that this requires additional simulations and that the affine parameter dependence cannot be accounted for, i.e., the matrix cannot be modified quickly to new parameters. 
\\

Equation \eqref{eq:FOF:qic} implies 
\begin{equation} \label{eq:App_simul:op_dft}
    \left(\td{I} -e^{(\td{A} - i \omega_k \td{I})\Delta t} \right)^{-1} \left(\td{I} - e^{\td{A}T} \right) \vd{v}_0 = \text{DFT}_k\left[ e^{\td{A}t_\set{J}} \vd{v}_0 \right] \text{.}
\end{equation}
That is, the effect of the operator on the left-hand side above applied to a vector $\vd{v}_0 \in \mathbb{C}^{N_x}$ is to take the $k$-th frequency of the DFT of the time series generated by using $\vd{v}_0$ as the initial condition to the linear system $\dot{\vd{v}} = \td{A} \td{v}$. 
\\

Using this equivalence, the $l$-th column of the matrix $\SPOD{k}^* \td{W} \left(\td{I} -e^{(\td{A} - i \omega_k \td{I})\Delta t} \right)^{-1} \left(\td{I} - e^{\td{A}T} \right) \td{\Phi}$ is $\SPOD{k}^* \td{W} \vdh{\phi}_k^l$, where $\vdh{\phi}_k^l = \text{DFT}_k[ e^{\td{A}t_\set{J}}  \vd{\phi}^l]$ is the $k$-th component of the DFT of the time series created by using the $l$-th column of $\td{\Phi}$ as an initial condition to the unforced linear system. Thus, by initializing this linear system with each columns of $\td{\Phi}$, and taking the DFT of each result, all matrices $\SPOD{k}^* \td{W} \left(\td{I} -e^{(\td{A} - i \omega_k \td{I})\Delta t} \right)^{-1} \left(\td{I} - e^{\td{A}T} \right) \td{\Phi}$ can be computed up to the time-stepping errors in integrating the linear system.
\\

In principle, linear runs could also be used to compute the action of $\SPOD{k}^* \td{W} \td{R}_k$, and in a different context, such a strategy has been employed by Ref. \cite{Farghadan25}. In the present context, such an approximation would involve performing linear runs with harmonic forcings with spatial modes corresponding to the vectors to which $\SPOD{k}^* \td{W} \td{R}_k$ is applied. However, unlike the runs described above, these linear runs would need to be long enough for transients to decay \cite{Farghadan25}. Moreover, with the triadic-interaction handling of the nonlinearity, there are prohibitively many vectors that would need to be used as the forcing in separate runs. Thus, we do not recommend this approach in this context.

\section{Algorithms} \label{app:algorithms}
The offline phase of the method is given in Algorithm~\ref{alg:offline}. The first step in the online phase of the method is computing the constant term using Algorithm~\ref{alg:online_const}. Then, either the fixed point iteration or the pseudo-time-stepping method where the nonlinear term $\vd{w}_\set{K}(\vd{a}_\set{K})$ is calculated using either Algorithm~\ref{alg:online_nonquad} or Algorithm~\ref{alg:online_quad}, depending on the nature of the nonlinearity.  For readability, we have written Algorithm~\ref{alg:online_quad} with for loops, but we note that in our numerical implementation, we used sparse matrices, which will run faster on most systems. 
\\

\begin{algorithm} 
\caption{SSOP (offline)}\label{alg:offline}
\begin{algorithmic}[1]
\State $\td{\Phi} = \texttt{POD}(\td{Q}^t,\td{W},p_1)$ \Comment{Intermediary basis} 
\State $\tdh{Q}_\set{K} = \texttt{WelchBlocks}(\td{Q}^t,N_\omega)$ \Comment{Obtain data matrices of each frequency from snapshot matrix}
\If{$\vd{n}$ \ \text{is non-quadratic}}
\State $[\td{U}^n, \td{P}^{nT}] = \texttt{DEIM}(\vd{n}(\td{Q}^t),p_2)$ \Comment{Obtain sample points and basis for nonlinearity}
\EndIf
\State $\left[ \SPOD{k}^{N_d},\td{\Lambda}_k^{N_d},\SPOD{k}^{},\td{\Lambda}_k^{} \right] = \texttt{SPOD}(\tdh{Q}_\set{K} ,r) $ \Comment{Get SPOD modes and energies}
\For{$k \in \set{K}$}
\State $\td{L}_k \gets \left(i\omega_k \td{I} - \td{A} \right)$ 
\State $\td{G}_k \gets \td{L}_k \tdh{Q}_k$  \Comment{Used for approximating resolvent}
\State $\td{E}_k \gets \SPOD{k}^{*} \td{W} \tdh{Q}_k (\td{G}_k^* \td{W} \td{G}_k)^{-1} \td{G}_k^* \td{W} \td{B}$ \Comment{Output variable}
\State $\td{J}_k \gets \td{\Phi}^* \td{W} \tdh{Q}_k (\td{G}_k^* \td{W} \td{G}_k)^{-1} \td{G}_k^* \td{W} \td{B}$ \Comment{Output variable}
\State $\tilde{\td{A}}_k \gets \SPOD{k}^{N_d*} \td{W} \td{A}\SPOD{k}^{N_d}$ \Comment{Used for calculating $\td{H}_k$} 
\State $\td{P}_k \gets \left[ \td{I}_{r_k \times r_k} \ \td{0}_{r_k \times (N_d - r_k)} \right]$ \Comment{Used for calculating $\td{H}_k$}
\State $\td{H}_k \gets \td{P}_k \left( \td{I} - \exp \left[ \left( \tilde{\td{A}} - i\omega_k \td{I} \right) \Delta t \right] \right)^{-1} \left( \td{I} - \exp \left[ \tilde{\td{A}}T \right] \right) \left( \SPOD{k}^{N_d*}\td{W}\td{\Phi} \right)$ \Comment{Output variable}

\If{$\vd{n}$ \ \text{is non-quadratic}}
\State $\td{N}_k \gets \SPOD{k}^{*} \td{W} \tdh{Q}_k (\td{G}_k^* \td{W} \td{G}_k)^{-1} \td{G}_k^* \td{W} \td{U}^n \left( \td{P}^{nT}\td{U}^n \right)^{-1}$ \Comment{Output variable}
\State $\td{M}_k \gets \td{\Phi}^* \td{W} \tdh{Q}_k (\td{G}_k^* \td{W} \td{G}_k)^{-1} \td{G}_k^* \td{W} \td{U}^n \left( \td{P}^{nT}\td{U}^n \right)^{-1}$ \Comment{Output variable}
\State $\td{S}_k \gets \td{P}^{nT} \SPOD{k}$ \Comment{Output variable}
\Else
\State $\vd{b}(\vd{q}_1,\vd{q}_2) := \frac{1}{2}\left[ \vd{n}(\vd{q}_1 + \vd{q}_2) - \vd{n}(\vd{q}_1) - \vd{n}(\vd{q}_2)   \right] $ \Comment{Symmetric bilinear form}
\State $\set{N}_k \gets \{ \}$ 
\For{$(l,i) \text{ such that } \omega_l + \omega_i = \omega_k$} \Comment{Looping over triadic interactions}
\For {$m \in \{ 1, \dots , r_l \}$}
\For {$n \in \{ 1, \dots , r_i \}$}
\State $\tilde{\vd{n}} \gets  \SPOD{k}^{*} \td{W} \tdh{Q}_k (\td{G}_k^* \td{W} \td{G}_k)^{-1} \td{G}^* \td{W} \vd{b}(\vd{\psi}_{lm},\vd{\psi}_{in}) $
\If{$\|\tilde{\vd{n}} \|_{x}^2 \lambda_{lm} \lambda_{in} > \epsilon$} \Comment{Checking if interaction meets threshold}
\State $\vd{n}_{klmn} \gets \tilde{\vd{n}}$ \Comment{Output variable}
\State $\vd{m}_{klmn} \gets \td{\Phi}^* \td{W} \tdh{Q}_k (\td{G}_k^* \td{W} \td{G}_k)^{-1} \td{G}^* \td{W} \vd{b}(\vd{\psi}_{lm},\vd{\psi}_{in})$ \Comment{Output variable}
\State $\set{N}_k \gets \set{N}_k \cup (l,m,n)$ \Comment{Updating list of interactions} 
\EndIf
\EndFor
\EndFor
\EndFor
\EndIf
\EndFor
\end{algorithmic}
\vspace{1\baselineskip}
\textbf{Inputs:} $\td{Q}^t$, time series data; $\Delta t$, time step used in data; $N_\omega$, number of time steps desired for solution blocks; $\vd{n}$, nonlinear function; $p_1$, number of modes to use in intermediary basis; $p_2$ (if DEIM is used), number of sample points to use in DEIM approximation; $r$, average number of modes per frequency for ROM; $\td{A}$, $\td{B}$, system matrices; $\td{W}$, weight matrix; $\epsilon$ (if triadic interactions are used), tolerance for inclusion.\\
\textbf{Outputs: } $\td{\Phi}$, intermediary basis; $\td{J}_\set{K}$, $\td{J}$ operators at all frequencies; $\td{E}_\set{K}$, $\td{E}$ operators at all frequencies ; $\td{H}_\set{K}$, $\td{H}$ operators at all frequencies. If using DEIM,  $\td{S}_\set{K}$, the sampling operators at all frequencies; $\td{N}_\set{K}$, the first set of hyper-reduction matrices; $\td{M}_\set{K}$ the second set of hyper-reduction matrices. Otherwise, $\set{N}_\set{K}$, the set $\set{N}$ at all frequencies; $\vd{n}_{klmn}$ for all $k \in \set{K}$, $(l,m,n) \in \set{N}_k$; $\vd{m}_{klmn}$ for all $k \in \set{K}$, $(l,m,n) \in \set{N}_k$.
\end{algorithm}

\begin{algorithm} 
\caption{SSOP (online, constant term)}\label{alg:online_const}
\begin{algorithmic}[1]
\State $\vdh{f}_\set{K} = \texttt{FFT}_\set{K}[\vd{f}_\set{J}]$ \Comment{FFT of forcing}
\State $\vd{q}_0^{\td{\Phi}} \gets \td{\Phi}^*\td{W} \vd{q}_0$ \Comment{The initial condition in the intermediary basis}
\State $\vdt{q}_0^{\td{\Phi}} \gets \vd{0}_{p_1 \times 1}$ \Comment{Initializing the forcing sum term}
\For{$k \in \set{K}$}
\State $\vdt{q}_0^{\td{\Phi}} \gets \vdt{q}_0^{\td{\Phi}} + \frac{1}{N_\omega}\td{J}_k \vdh{f}_k$ \Comment{Influence of each frequency on forcing sum}
\EndFor

\For{$k \in \set{K}$}
\State $\vd{c}_k \gets \td{E}_k \vdh{f}_k + \td{H}_k\left( \vd{q}_0^{\td{\Phi}} -  \vdt{q}_0^{\td{\Phi}} \right) $ \Comment{Constructing each frequency of the constant term}
\EndFor

\end{algorithmic}
\vspace{1\baselineskip}
\textbf{Inputs:} $\vd{q}_0$, the initial condition; $\td{\Phi}$, intermediary basis (POD modes); $\vd{W}$, weight matrix; $\vd{f}_\set{J}$, the forcing on the temporal grid; $\td{J}_\set{K}$, $\td{J}$ operators at all frequencies (defined in \ref{eq:Lin:Jdef}); $\td{E}_\set{K}$, $\td{E}$ operators at all frequencies (defined in \eqref{eq:Lin:Edef}); $\td{H}_\set{K}$, $\td{H}$ operators at all frequencies (defined below \ref{eq:Lin:ROMfinal}).\\
\textbf{Outputs:} $\vd{c}_\set{K}$, all frequency components of the constant term.
\end{algorithm}

\begin{algorithm} 
\caption{SSOP (online, non-quadratic nonlinear term)}\label{alg:online_nonquad}
\begin{algorithmic}[1]
\For{$k \in \set{K}$}
\State $\vdt{q}^s_k \gets \td{S}_k \vdt{a}_k$ \Comment{Getting sample points in frequency domain}
\EndFor
\State $\vd{q}_\set{J}^s = \texttt{IFFT}[ \vdt{q}_\set{K}^s ]$ \Comment{Sample points in time domain}
\State $\vd{n}_\set{J}^s = \vd{n}(\vd{q}_\set{J}^s)$ \Comment{Nonlinearity at sample points in time domain}
\State $\vdh{n}_\set{K}^s = \texttt{FFT} [ \vd{n}_\set{J}^s ] $ \Comment{Nonlinearity at sample points in frequency domain}

\State $\vdt{q}_0^{\td{\Phi}} \gets \vd{0}_{p_1 \times 1}$
\For{$k \in \set{K}$}
\State $\vdt{q}_0^{\td{\Phi}} \gets \vdt{q}_0^{\td{\Phi}} + \frac{1}{N_\omega} \td{M}_k \vdh{n}_k^s$
\State $\vd{w}_k \gets \td{N}_k \vdh{n}_k^s$
\EndFor

\For{$k \in \set{K}$}
\State $\vd{w}_{k} \gets \vd{w}_{k} - \td{H}_k \vdt{q}_0^{\td{\Phi}}$
\EndFor

\end{algorithmic}
\vspace{1\baselineskip}
\textbf{Inputs:} $\vd{a}_\set{K}$, the SPOD coefficients at the current iteration; $\td{S}_\set{K}$, the sampling operators at all frequencies; $\vd{n}$, the nonlinear function; $\td{N}_\set{K}$, the first set of matrices hyper-reduction matrices (see \eqref{eq:nonlinear:deim:matrix1}); $\td{M}_\set{K}$ the second set of hyper-reduction matrices (see \eqref{eq:nonlinear:deim:matrix2}); $\td{H}_\set{K}$, the set of matrices defined under \eqref{eq:Lin:ROMfinal}.\\
\textbf{Output:} $\vd{w}_\set{K}$, the approximate nonlinear term formed by concatenating $\vd{w}$ at all frequencies.
\end{algorithm}

\begin{algorithm} 
\caption{SSOP (online, quadratic nonlinear term)}\label{alg:online_quad}
\begin{algorithmic}[1]

\State $\vd{w}^{\td{\Phi}} \gets \vd{0}_{p_1 \times 1}$ 
\For{$k \in \set{K}$} 
\State $\vd{w}_k \gets \vd{0}_{r_k \times 1}$
\For{$(l,m,n) \in \set{N}_k$}
\State $\vd{w}_k \gets \vd{w}_k + \vd{n}_{klmn} \tilde{a}_{lm}\tilde{a}_{in}$
\State $\vd{w}^{\td{\Phi}} \gets \vd{w}^{\td{\Phi}} + \vd{m}_{klmn} \tilde{a}_{lm}\tilde{a}_{in}$
\EndFor
\EndFor

\For{$k \in \set{K}$}
\State $\vd{w}_k \gets  \vd{w}_k - \td{H}_k \vd{w}^{\td{\Phi}}$
\EndFor

\end{algorithmic}
\vspace{1\baselineskip}
\textbf{Inputs:} $\vd{a}_\set{K}$, the SPOD coefficients at the current iteration; $\set{N}_\set{K}$, the set $\set{N}$ at all frequencies (see \eqref{eq:nonlinear:Ndef}); $\vd{n}_{klmn}$ for all $k \in \set{K}$, $(l,m,n) \in \set{N}_k$ (see \eqref{eq:nonlinear:nklmndef}); $\vd{m}_{klmn}$ for all $k \in \set{K}$, $(l,m,n) \in \set{N}_k$ (see \eqref{eq:nonlinear:mklmndef});  $\td{H}_\set{K}$, the set of matrices defined under \eqref{eq:Lin:ROMfinal}.\\
\textbf{Output:} $\vd{w}_\set{K}$, the approximate nonlinear term formed by concatenating $\vd{w}$ at all frequencies.

\end{algorithm}

\end{appendices}
\pagebreak

\bibliographystyle{abbrv}
\bibliography{biblio}

\begin{thebibliography}{10}

\bibitem{Anderson65}
D.~G. Anderson.
\newblock Iterative procedures for nonlinear integral equations.
\newblock {\em Journal of the ACM}, 12(4):547–560, Oct. 1965.

\bibitem{Aubry88}
N.~Aubry, P.~Holmes, J.~L. Lumley, and E.~Stone.
\newblock The dynamics of coherent structures in the wall region of a turbulent boundary layer.
\newblock {\em Journal of Fluid Mechanics}, 192:115–173, July 1988.

\bibitem{Bagheri09}
S.~Bagheri, D.~S. Henningson, J.~H{\oe}pffner, and P.~J. Schmid.
\newblock Input-output analysis and control design applied to a linear model of spatially developing flows.
\newblock {\em Applied Mechanics Reviews}, 62(2), Feb. 2009.

\bibitem{Brunton14}
S.~L. Brunton, J.~H. Tu, I.~Bright, and J.~N. Kutz.
\newblock Compressive sensing and low-rank libraries for classification of bifurcation regimes in nonlinear dynamical systems.
\newblock {\em SIAM Journal on Applied Dynamical Systems}, 13(4):1716–1732, Jan. 2014.

\bibitem{Carlberg10}
K.~Carlberg, C.~Bou‐Mosleh, and C.~Farhat.
\newblock Efficient non‐linear model reduction via a least‐squares {P}etrov–{G}alerkin projection and compressive tensor approximations.
\newblock {\em International Journal for Numerical Methods in Engineering}, 86(2):155–181, Oct. 2010.

\bibitem{Chaturantabut10}
S.~Chaturantabut and D.~C. Sorensen.
\newblock Nonlinear model reduction via discrete empirical interpolation.
\newblock {\em SIAM Journal on Scientific Computing}, 32(5):2737–2764, Jan. 2010.

\bibitem{Chen11}
K.~K. Chen and C.~W. Rowley.
\newblock Optimal actuator and sensor placement in the linearised complex {G}inzburg{\textendash}{L}andau system.
\newblock {\em Journal of Fluid Mechanics}, 681:241--260, June 2011.

\bibitem{Chen21}
Y.~Chen and C.~Liu.
\newblock Data driven robust control of complex ginzburg-landau equations for spatial developing flow.
\newblock In {\em 2021 IEEE 4th International Conference on Automation, Electronics and Electrical Engineering (AUTEEE)}, page 490–500. IEEE, Nov. 2021.

\bibitem{Choi21}
Y.~Choi, P.~Brown, W.~Arrighi, R.~Anderson, and K.~Huynh.
\newblock Space–time reduced order model for large-scale linear dynamical systems with application to {B}oltzmann transport problems.
\newblock {\em Journal of Computational Physics}, 424:109845, 2021.

\bibitem{Choi19}
Y.~Choi and K.~Carlberg.
\newblock Space--time least-squares petrov--galerkin projection for nonlinear model reduction.
\newblock {\em SIAM Journal on Scientific Computing}, 41(1):A26--A58, 2019.

\bibitem{Farghadan25}
A.~Farghadan, E.~Martini, and A.~Towne.
\newblock Scalable resolvent analysis for three-dimensional flows.
\newblock {\em Journal of Computational Physics}, 524:113695, Mar. 2025.

\bibitem{Frame25}
P.~Frame, C.~Lin, O.~T. Schmidt, and A.~Towne.
\newblock Linear model reduction using spectral proper orthogonal decomposition.
\newblock {\em Computer Methods in Applied Mechanics and Engineering}, 447:196–215, Dec. 2025.

\bibitem{Frame23}
P.~Frame and A.~Towne.
\newblock Space-time {POD} and the {H}ankel matrix.
\newblock {\em {PLOS} {ONE}}, 18(8):e0289637, Aug. 2023.

\bibitem{Hall13}
K.~C. Hall, K.~Ekici, J.~P. Thomas, and E.~H. Dowell.
\newblock Harmonic balance methods applied to computational fluid dynamics problems.
\newblock {\em International Journal of Computational Fluid Dynamics}, 27(2):52–67, Jan. 2013.

\bibitem{Hall02}
K.~C. Hall, J.~P. Thomas, and W.~S. Clark.
\newblock Computation of unsteady nonlinear flows in cascades using a harmonic balance technique.
\newblock {\em AIAA Journal}, 40(5):879–886, May 2002.

\bibitem{Hoang22}
C.~Hoang, K.~Chowdhary, K.~Lee, and J.~Ray.
\newblock Projection-based model reduction of dynamical systems using space–time subspace and machine learning.
\newblock {\em Computer Methods in Applied Mechanics and Engineering}, 389:114341, Feb. 2022.

\bibitem{Ilak10}
M.~Ilak, S.~Bagheri, L.~Brandt, C.~W. Rowley, and D.~S. Henningson.
\newblock Model reduction of the nonlinear complex ginzburg–landau equation.
\newblock {\em SIAM Journal on Applied Dynamical Systems}, 9(4):1284–1302, Jan. 2010.

\bibitem{Choi21b}
Y.~Kim, K.~Wang, and Y.~Choi.
\newblock Efficient space–time reduced order model for linear dynamical systems in python using less than 120 lines of code.
\newblock {\em Mathematics}, 9(14), 2021.

\bibitem{Kramer24}
B.~Kramer, B.~Peherstorfer, and K.~E. Willcox.
\newblock Learning nonlinear reduced models from data with operator inference.
\newblock {\em Annual Review of Fluid Mechanics}, 56(1):521–548, Jan. 2024.

\bibitem{Li25}
X.~Li and D.~Lasagna.
\newblock Space-time nonlinear reduced-order modelling for unsteady flows, 2025.

\bibitem{Lin19}
C.~Lin.
\newblock Model order reduction in the frequency domain via spectral proper orthogonal decomposition.
\newblock {\em Master's thesis, University of Illinois at Urbana-Champaign}, 2019.

\bibitem{Lumley67}
J.~L. Lumley.
\newblock The structure of inhomogeneous turbulent flows.
\newblock {\em Atmospheric Turbulence and Radio Wave Propagation}, 1967.

\bibitem{Lumley70}
J.~L. Lumley.
\newblock Stochastic tools in turbulence.
\newblock 1970.

\bibitem{McKeon17}
B.~J. McKeon.
\newblock The engine behind (wall) turbulence: perspectives on scale interactions.
\newblock {\em Journal of Fluid Mechanics}, 817:P1, 2017.

\bibitem{McKeon10}
B.~J. McKeon and A.~S. Sharma.
\newblock A critical-layer framework for turbulent pipe flow.
\newblock {\em Journal of Fluid Mechanics}, 658:336--382, July 2010.

\bibitem{Moore81}
B.~Moore.
\newblock Principal component analysis in linear systems: Controllability, observability, and model reduction.
\newblock {\em IEEE Transactions on Automatic Control}, 26(1):17–32, Feb. 1981.

\bibitem{Noack03}
B.~R. Noack, K.~Afanasiev, M.~Morzyński, G.~Tadmor, and F.~Thiele.
\newblock A hierarchy of low-dimensional models for the transient and post-transient cylinder wake.
\newblock {\em Journal of Fluid Mechanics}, 497:335–363, Dec. 2003.

\bibitem{Oja82}
E.~Oja.
\newblock Simplified neuron model as a principal component analyzer.
\newblock {\em Journal of Mathematical Biology}, 15(3):267–273, Nov. 1982.

\bibitem{Orszag71}
S.~A. Orszag.
\newblock On the elimination of aliasing in finite-difference schemes by filtering high-wavenumber components.
\newblock {\em Journal of the Atmospheric Sciences}, 28(6):1074–1074, Sept. 1971.

\bibitem{Padovan24}
A.~Padovan, B.~Vollmer, and D.~J. Bodony.
\newblock Data-driven model reduction via non-intrusive optimization of projection operators and reduced-order dynamics, 2024.

\bibitem{Parish21}
E.~J. Parish and K.~T. Carlberg.
\newblock Windowed least-squares model reduction for dynamical systems.
\newblock {\em Journal of Computational Physics}, 426:109939, 2021.

\bibitem{Peherstorfer16}
B.~Peherstorfer and K.~Willcox.
\newblock Data-driven operator inference for nonintrusive projection-based model reduction.
\newblock {\em Computer Methods in Applied Mechanics and Engineering}, 306:196–215, July 2016.

\bibitem{rigas21}
G.~Rigas, D.~Sipp, and T.~Colonius.
\newblock Nonlinear input/output analysis: application to boundary layer transition.
\newblock {\em Journal of Fluid Mechanics}, 911, 2021.

\bibitem{Rowley05}
C.~W. Rowley.
\newblock Model reduction for fluids, using balanced proper orthogonal decomposition.
\newblock {\em International Journal of Bifurcation and Chaos}, 15(03):997--1013, Mar. 2005.

\bibitem{Schmid07}
P.~J. Schmid.
\newblock Nonmodal stability theory.
\newblock {\em Annual Review of Fluid Mechanics}, 39(1):129–162, Jan. 2007.

\bibitem{Schmidt20}
O.~T. Schmidt.
\newblock Bispectral mode decomposition of nonlinear flows.
\newblock {\em Nonlinear Dynamics}, 102(4):2479–2501, Nov. 2020.

\bibitem{Sipp20}
D.~Sipp, M.~Fosas~de Pando, and P.~J. Schmid.
\newblock Nonlinear model reduction: A comparison between pod-galerkin and pod-deim methods.
\newblock {\em Computers \& Fluids}, 208:104628, Aug. 2020.

\bibitem{Sirovich87}
L.~Sirovich.
\newblock Turbulence and the dynamics of coherent structures. ii. symmetries and transformations.
\newblock {\em Quarterly of Applied Mathematics}, 45(3):573–582, 1987.

\bibitem{Towne21}
A.~Towne.
\newblock Space-time {G}alerkin projection via spectral proper orthogonal decomposition and resolvent modes.
\newblock In {\em AIAA Scitech 2021 Forum}. American Institute of Aeronautics and Astronautics, Jan. 2021.

\bibitem{Towne2018spectral}
A.~Towne, O.~T. Schmidt, and T.~Colonius.
\newblock Spectral proper orthogonal decomposition and its relationship to dynamic mode decomposition and resolvent analysis.
\newblock {\em Journal of Fluid Mechanics}, 847:821--867, 2018.

\bibitem{Trefethen93}
L.~N. Trefethen, A.~E. Trefethen, S.~C. Reddy, and T.~A. Driscoll.
\newblock Hydrodynamic stability without eigenvalues.
\newblock {\em Science}, 261(5121):578--584, 1993.

\bibitem{Walker11}
H.~F. Walker and P.~Ni.
\newblock Anderson acceleration for fixed-point iterations.
\newblock {\em SIAM Journal on Numerical Analysis}, 49(4):1715–1735, Jan. 2011.

\bibitem{Welch67}
P.~Welch.
\newblock The use of fast fourier transform for the estimation of power spectra: A method based on time averaging over short, modified periodograms.
\newblock {\em IEEE Transactions on Audio and Electroacoustics}, 15(2):70--73, 1967.

\bibitem{Willcox02}
K.~Willcox and J.~Peraire.
\newblock Balanced model reduction via the proper orthogonal decomposition.
\newblock {\em {AIAA} Journal}, 40(11):2323--2330, Nov. 2002.

\bibitem{Yano14}
M.~Yano, A.~T. Patera, and K.~Urban.
\newblock A space-time hp-interpolation-based certified reduced basis method for {B}urgers{\textquotesingle} equation.
\newblock {\em Mathematical Models and Methods in Applied Sciences}, 24(09):1903--1935, May 2014.

\end{thebibliography}

\end{document}